\documentclass[titlepage,draft,12pt]{article} 
\usepackage{amsfonts,latexsym} 
\usepackage{pstricks,pst-plot}
\usepackage[a4paper]{geometry}

\date{}

\newcommand{\ep}{\varepsilon}
\newcommand{\qed}{{\penalty 10000\mbox{$\quad\Box$}}}
\newcommand{\re}{\mathbb{R}}

\newcommand{\n}{\mathbb{N}}

\newcommand{\ul}{u_{\lambda}}
\newcommand{\vl}{v_{\lambda}}
\newcommand{\vle}{v_{\lambda,\varepsilon}}
\newcommand{\al}{a_{\lambda}}
\newcommand{\bl}{b_{\lambda}}
\newcommand{\fl}{f_{\lambda}}
\newcommand{\tl}{\tau_{\lambda}}
\newcommand{\phil}{\varphi_{\lambda}}
\newcommand{\phile}{\varphi_{\lambda,\varepsilon}}
\newcommand{\psie}{\psi_{\varepsilon}}
\newcommand{\lk}{\lambda_{k}}
\newcommand{\laxi}{\lambda(\xi)}

\newcommand{\zt}{\widehat{z}}
\newcommand{\ut}{\widehat{u}}
\newcommand{\vt}{\widehat{v}}
\newcommand{\wt}{\widehat{w}}
\newcommand{\ft}{\widehat{f}}
\newcommand{\vft}{\widehat{\varphi}}

\newtheorem{thm}{Theorem}[section]

\newtheorem{rmk}[thm]{Remark}

\newtheorem{lemma}[thm]{Lemma}

\title{Local and global smoothing effects for some linear hyperbolic 
equations with a strong dissipation}

\author{Marina Ghisi\vspace{1ex}\\ 
{\normalsize Universit\`a degli Studi di Pisa} \\
{\normalsize Dipartimento di Matematica}\\ 
{\normalsize PISA (Italy)}\\
{\normalsize e-mail: \texttt{ghisi@dm.unipi.it}}
\and
Massimo Gobbino\vspace{1ex}\\ 
{\normalsize Universit\`a degli Studi di Pisa} \\
{\normalsize Dipartimento di Matematica}\\ 
{\normalsize PISA (Italy)}\\  
{\normalsize e-mail: \texttt{m.gobbino@dma.unipi.it}}
\and
Alain Haraux\vspace{1ex}\\ 
{\normalsize Universit\'{e} Pierre et Marie Curie} \\
{\normalsize Laboratoire Jacques-Louis Lions}\\ 
{\normalsize PARIS (France)}\\  
{\normalsize e-mail: \texttt{haraux@ann.jussieu.fr}}}

\begin{document}
\maketitle
\begin{abstract}
	We consider an abstract second order linear equation with a 
	strong dissipation, namely a friction term which depends on a 
	power of the ``elastic'' operator.
	
	In the homogeneous case, we investigate the phase spaces in which 
	the initial value problem gives rise to a semigroup, and the 
	further regularity of solutions. In the non-homogeneous case, we 
	study how the regularity of solutions depends on the 
	regularity of forcing terms, and we characterize the spaces 
	where a bounded forcing term yields a bounded solution.
	
	What we discover is a variety of different regimes, with completely 
	different behaviors, depending on the exponent in the friction 
	term.
	
	We also provide counterexamples in order to show the optimality 
	of our results.
	
\vspace{1cm}

\noindent{\textbf Mathematics Subject Classification 2010 (MSC2010):}
35L10, 35L15, 35L20.

% 35L10: Second-order hyperbolic equations
% 35L15: Initial value problems for second-order hyperbolic equations
% 35L20: Initial-boundary value problems for second-order hyperbolic equations

\vspace{1cm} 

\end{abstract}

%%%%%%%%%%%%%%%%%%%%%
%                   %
%   Inizio lavoro   %
%                   %
%%%%%%%%%%%%%%%%%%%%%
 
\section{Introduction}

Let $H$ be a separable real Hilbert space.  For every $x$ and $y$ in
$H$, $|x|$ denotes the norm of $x$, and $\langle x,y\rangle$ denotes
the scalar product of $x$ and $y$.  Let $A$ be a self-adjoint linear
operator on $H$ with dense domain $D(A)$.  We assume that $A$ is
nonnegative, namely $\langle Ax,x\rangle\geq 0$ for every $x\in D(A)$,
so that for every $\alpha\geq 0$ the power $A^{\alpha}x$ is defined
provided that $x$ lies in a suitable domain $D(A^{\alpha})$.

We consider the second order linear evolution equation
\begin{equation}
	u''(t)+2\delta A^{\sigma}u'(t)+Au(t)=f(t),
	\quad\quad
	t\geq 0,
	\label{pbm:eqn}
\end{equation}
where $\delta>0$, $\sigma\geq 0$, and $f:[0,+\infty)\to H$, with 
initial data
\begin{equation}
	u(0)=u_{0},
	\hspace{3em}
	u'(0)=u_{1}.
	\label{pbm:data}
\end{equation}

Several wave equations fit in this abstract framework, for example
\begin{eqnarray}
	 & u_{tt}+u_{t}-\Delta u= f(t, x)
	 \hspace{3em} 
	 (A=-\Delta,\quad \sigma=0), & 
	\nonumber  \\
	\noalign{\vspace{1ex}}
	 & u_{tt}-\Delta u_{t}-\Delta u= f(t, x)
	 \hspace{3em} 
	 (A=-\Delta,\quad\sigma=1), & 
	\label{concrete:v-e}  \\
	\noalign{\vspace{1ex}}
	 & u_{tt}-\Delta u_{t}+\Delta^{2} u=f(t, x)
	 \hspace{3em} 
	 (A=\Delta^{2},\quad\sigma=1/2), & 
	 \nonumber
\end{eqnarray}
with ad hoc boundary conditions.  The case $\sigma=0$ is the standard
damped wave equation, the case $\sigma=1$ is sometimes called
visco-elastic damping, the case $\sigma\in(0,1)$ is usually referred
to as structural damping or fractional damping.  The case $\sigma>1$
seems to be quite unexplored.

Mathematical models of this kind were proposed in~\cite{CR}, and then
rigorously analyzed by many authors from different points of view.  In
the abstract setting, a natural idea is to set $U(t):=(u(t),u'(t))$,
so that one can interpret (\ref{pbm:eqn}) as a first order system
$$U'(t)+\mathcal{A}_{\sigma}U(t)=\mathcal{F}(t),$$
where
$$\mathcal{A}_{\sigma}=\left(
\begin{array}{cc}
	0 & -I  \\
	A & 2\delta A^{\sigma}
\end{array}
\right)$$
acts on some product space $\mathcal{H}$, usually chosen equal to
$D(A^{1/2})\times H$, and 
$$\mathcal{F}(t)=\left(
\begin{array}{c}
	0  \\
	f(t)
\end{array}
\right).$$  

In this framework, several papers have been devoted to properties of
semigroup generated by the operator $\mathcal{A}_{\sigma}$, such as
analyticity (true in the case $\sigma \geq 1/2$) or Gevrey regularity
(true in the case $0<\sigma<1/2$) .  The interested reader is referred
to~\cite{CT1,CT2,CT3}, and to the more recent
papers~\cite{FGR,HO,LL,mugnolo}.  In these papers $A$ is usually
assumed to be strictly positive, $\sigma\in[0,1]$, and the phase space
is $D(A^{1/2})\times H$.

On a completely different side, the community working on dispersive
hyperbolic equations considered the concrete equation
(\ref{concrete:v-e}) (with $f=0$), or its $\sigma$-generalization, in
the whole space $\re^{n}$ or in suitable classes of unbounded domains,
obtaining $L^{p}$-$L^{q}$ estimates or energy decay estimates.  The
interested reader is referred to~\cite{shibata,I1,I2,I3} and to the
references quoted therein.  In these papers $A$ is just nonnegative,
$\sigma\in[0,1]$, and the phase space is once again $D(A^{1/2})\times
H$, with some $L^{p}$ restrictions, but the dispersive properties of
$A$ are essential.

We conclude this brief historical survey, which is far from being
complete, by mentioning the existence of some literature on nonlinear
wave equations with fractional damping (see for
example~\cite{nishihara,ono-nishihara} or the more recent
work~\cite{DAR}).

In this paper we consider the abstract problem
(\ref{pbm:eqn})--(\ref{pbm:data}) in its full generality, with the aim
of providing a complete picture in the whole range $\sigma\geq 0$.

We begin our study by considering the homogeneous case where
$f(t)\equiv 0$.  The first question we address is the choice of the
phase space.  Having in mind the standard setting for the
non-dissipative case ($\delta=0$), one is naturally led to consider
the phase space $D(A^{1/2})\times H$, or more generally
$D(A^{\alpha})\times D(A^{\alpha-1/2})$, namely with ``gap $1/2$''.
In the non-dissipative case, this choice has physical motivations (it
is the usual ``energy space''), but it is also dictated by the
equation itself in the sense that the initial value problem generates
a continuous semigroup in the phase space $D(A^{\alpha_{0}})\times
D(A^{\alpha_{1}})$ if and only if $\alpha_{0}-\alpha_{1}=1/2$.  We say
that 1/2 is the ``phase space gap''.

Also further time-derivatives of $u$ respect this gap, in the sense
that the $m$-th time-derivative $u^{(m)}(t)$ lies in the space
$D(A^{\alpha-(m-1)/2})$ for every $m\geq 1$ such that $\alpha\geq
(m-1)/2$, and once again the exponents are optimal if $A$ is
unbounded.  Thus we say that 1/2 is also the ``derivative gap''.  We
stress that in both cases the value 1/2 is chosen by the equation
itself.  \medskip

What about the dissipative case?  In statements~(1)
and~(2) of Theorem~\ref{thm:homog} we investigate phase space gaps and
derivative gaps.  Two different regimes appear.
\begin{itemize}
	\item As long as $0\leq\sigma\leq 1/2$, both the phase space gap
	and the derivative gap are equal to 1/2.  Thus from this point of
	view the picture is exactly the same as in the non-dissipative
	case.

	\item For $\sigma>1/2$ things are different.  First of all, the
	initial value problem generates a continuous semigroup in
	$D(A^{\alpha_{0}})\times D(A^{\alpha_{1}})$ if and only if
	$1-\sigma\leq\alpha_{0}-\alpha_{1}\leq\sigma$, namely there is an
	interval of possible phase space gaps.  This interval is always
	centered in 1/2, and contains also negative values when
	$\sigma>1$.
	
	As for further time-derivatives, it turns out that $u^{(m)}(t)\in
	D(A^{\alpha_{1}-(m-1)\sigma})$ for every $m\geq 1$ such that
	$\alpha_{1}\geq (m-1)\sigma$, which means that in this regime we
	now have a derivative gap equal to $\sigma$.

\end{itemize}

The second question we address in the homogeneous case is the
regularity of solutions for $t>0$, since a strong dissipation is
expected to have a smoothing effect.  Here three regimes appear, as
shown by statements~(3) and~(4) of Theorem~\ref{thm:homog}.
\begin{itemize}
	\item For $\sigma=0$ there is no further regularity for $t>0$, as
	in the non-dissipative case.

	\item For $\sigma\in(0,1)$ there is an instantaneous smoothing
	effect similar to parabolic equations, in the sense that $u\in
	C^{\infty}((0,+\infty),D(A^{\alpha}))$ for every $\alpha\geq 0$.
	Nevertheless, from the quantitative point of view (usually stated
	in terms of Gevrey spaces), this effect is actually weaker than in
	the parabolic case.  We do not deepen this issue in the present
	paper, and we refer the interested reader to \cite{MR:gevrey}.

	\item For $\sigma\geq 1$ the dissipation is so strong that it
	prevents too much smoothing, but a new feature appears, namely
	$u^{(m)}\in C^{0}((0,+\infty),D(A^{\alpha_{0}+m(\sigma-1)}))$.
	Since $\sigma-1\geq 0$, this is the opposite of the classical
	regularity loss: the higher is the time-derivation order, the
	higher is the space regularity!  We stress that this is true only
	for positive times.  If we are interested in the regularity up to
	$t=0$, then even for $\sigma\geq 1$ there is regularity loss in
	the standard direction, with derivative gap equal to $\sigma$, as
	already observed.

\end{itemize}

After settling the homogeneous case, we study several properties in
presence of a non-trivial forcing $f(t)$.  In this case each solution
is the sum of the solution of the corresponding homogeneous equation
with the same initial data, and the solution of the non-homogeneous
equation with forcing $f(t)$ and null initial data.  Thus in the
non-homogenous case, by relying on the previous results, we are
reduced to study the special case $u_{0}=u_{1}=0$.  We address three
issues.

\begin{itemize}
	\item First of all, we consider a forcing term $f\in
	L^{\infty}((0,T),H)$ and we describe the spaces of the form $
	D(A^{\alpha})$ where the solution $u(t)$ and its derivative
	$u'(t)$ lie (see Theorem~\ref{thm:nh-local} and
	Remark~\ref{rmk:nh}).

	\item Then we consider a forcing term $f\in
	L^{\infty}((0,+\infty),H)$, defined and bounded for all positive
	times, and we characterize the spaces of the form $ D(A^{\alpha})$
	where $u(t)$ and $u'(t)$ are (globally) bounded.  The answer is
	given by Theorem~\ref{thm:nh-bound} (see also
	Remark~\ref{rmk:bounded}), and it is somewhat unexpected.  Indeed
	it turns out that $u'(t)$ is always globally bounded in all the
	spaces of the form $ D(A^{\alpha})$ to which it belongs, while
	$u(t)$ is globally bounded in all the spaces of the form $
	D(A^{\alpha})$ to which it belongs if and only if
	$\sigma\in[0,1]$.  On the contrary, when $\sigma>1$ we have that
	$u(t)\in D(A^{\alpha})$ for all $\alpha\leq\sigma$, but $u(t)$ is
	globally bounded in $D(A^{\alpha})$ only for $\alpha<1$.  

	\item As a third issue, we apply our techniques to a somewhat
	different question.  We consider the non-homogeneous
	equation~(\ref{pbm:eqn}) with a bounded forcing term $f\in
	L^{\infty}(\re,H)$ defined on the whole real line, and we ask
	ourselves whether there exists a solution which is globally
	bounded in some phase space.  

\end{itemize}

The third issue above is usually referred to as \emph{non-resonance
property}, and has been studied by many authors in the concrete case
of hyperbolic equations in a bounded domain $\Omega$ with linear or
nonlinear local dissipation terms (see for example \cite{AP, Bi, B-H,
H0, H4, H2, HZ} and the references therein).  Except when additional
conditions are assumed on $f$, such as anti-periodicity or more
regularity (see~\cite{Anti-P,Reg-Forc}), all these authors had to
assume, even in the case of a periodic forcing term, that the damping
operator carries $H^{1}_{0}(\Omega)$ to $H^{-1}(\Omega)$ (namely
$D(A^{1/2})$ to $D(A^{-1/2})$ in the abstract setting) in a bounded
manner.  What was not clear is whether this is a fundamental
obstruction or not.  In a different direction, a linear dissipation
term of the form $Bu'$ was considered in~\cite{h-brezil}, with the
assumption that $B$ carries $D(A^{1/2})$ to $H$ in a bounded manner.
The result proved in~\cite{h-brezil} is that one has the non-resonance
property if and only if all solutions of the corresponding homogeneous
equation decay to 0 in a uniform exponential way.  Once again, it was
not clear whether or not this uniform exponential decay is still a
fundamental requirement for more general dissipation terms.

After thinking it over for several decades without any clear
answer, it seems reasonable to consider the toy model where the 
damping is provided by a linear unbounded operator which does not 
carry $D(A^{1/2})$ to $D(A^{-1/2})$. This led us to equation 
(\ref{pbm:eqn}) with $\sigma>1$, which was actually the initial 
motivation of this paper. 

The answer is somewhat surprising.  Indeed in
Theorem~\ref{thm:glob-bound} we prove that (\ref{pbm:eqn}) has the
non-resonance property for every $\delta>0$ and every $\sigma\geq 0$.
For $\sigma>1$, this non-resonance result seems completely new and
makes a sharp contrast with the conclusions of~\cite{h-brezil} because
for $A$ unbounded and $\sigma> 1$ the semigroup is not exponentially
stable.  Moreover the bounded solution stays bounded in phase spaces
which are stronger than the usual energy spaces traditionally
considered in these problems, thus showing that linear overdamping
improves the non-resonance property.  A quite surprising phenomenon
which has to be better understood by looking at different damping
operators, linear and nonlinear.  
\medskip

Finally, all the exponents involved in our regularity and boundedness
results are optimal in general, as shown by the counterexamples of
Theorem~\ref{thm:optimal}.  As far as we know, no such counterexamples
were known before in the literature, and also most of the different
strategies used in the construction seem to be new.

\medskip

Most of the result obtained in this paper are based on a thorough
knowledge of the asymptotic behavior, as $\lambda\to +\infty$, of the
roots of the characteristic polynomial
\begin{equation}
	x^{2}+2\delta\lambda^{\sigma}x+\lambda.
	\label{char-pol}
\end{equation}

The form and asymptotic behavior of the roots is different in
different ranges of $\sigma$ (namely $\sigma=0$, $\sigma\in(0,1/2)$,
$\sigma=1/2$, $\sigma\in(1/2,1)$, $\sigma=1$, $\sigma>1$), giving rise
to the composite picture described in our results above.

%\bigskip 
We believe that these optimal regularity and boundedness
results might provide a benchmark when looking at more general
equations, for example nonlinear equations, or equations in which
$A^{\sigma}u'(t)$ is replaced by the more general friction term
$Bu'(t)$, where $B$ is an operator ``comparable'' with some power
$A^{\sigma}$ (as in the original models in~\cite{CR}).

\bigskip

This paper is organized as follows.  In Section~\ref{sec:results} we
state all our main results.  In Section~\ref{sec:proofs-autonomous} we
give the proofs for the homogeneous case $f= 0$ and in
Section~\ref{sec:non-autonomous} we give the proofs for the forced
case.  In Section~\ref{sec:counterexamples} we exhibit some examples
showing the optimality of our regularity and boundedness results in
the forced case.

\setcounter{equation}{0}
\section{The results}\label{sec:results}

Before stating our results, let us spend just a few words on the
notion of solution.  Weak solutions to evolution problems can be
introduced in several ways, for example through density arguments as
limits of classical solutions, or through integral forms or
distributional formulations.  All these notions are equivalent in the
case of a damped wave-type equation equation such as~(\ref{pbm:eqn}).

Moreover, thanks to the spectral theory for self-adjoint operators
(namely Fourier series or Fourier transform in concrete cases), the
study of (\ref{pbm:eqn})--(\ref{pbm:data}) reduces to the study of a
suitable family of ordinary differential equations (see
section~\ref{sec:notation} for further details).  In this way one can
prove the well-known results concerning existence of a unique solution
for quite general initial data and forcing terms (even distributions
or hyperfunctions), up to admitting that the solution takes its values
in a very large Hilbert space as well (once again distributions or
hyperfunctions).  Thus in the sequel we say ``the solution'' without
any further specification.

In this paper we investigate how the regularity of initial data and
forcing terms affects the regularity of solutions.  Let us start
with the homogeneous case.

\begin{thm}[The homogeneous equation]\label{thm:homog}
	Let $H$ be a separable Hilbert space, and let $A$ be a
	self-adjoint nonnegative operator on $H$ with dense domain $D(A)$.
	For every $\sigma\geq 0$, $\delta>0$, $\alpha_{0}\geq 0$,
	$\alpha_{1}\geq 0$, we consider the unique solution to the
	homogeneous linear equation
	\begin{equation}
		u''(t)+2\delta A^{\sigma}u'(t)+Au(t)=0,
		\quad\quad
		t\geq 0,
		\label{pbm:h-eqn}
	\end{equation}
	with initial data
	\begin{equation}
		u(0)=u_{0}\in D(A^{\alpha_{0}}),
		\quad\quad
		u'(0)=u_{1}\in D(A^{\alpha_{1}}).
		\label{pbm:h-data}
	\end{equation}
	
	Let us set $\gamma:=\max\{1/2,\sigma\}$, and let us assume that
	\begin{equation} 
		1-\gamma\leq \alpha_{0}-\alpha_{1}\leq
		\gamma.
		\label{gap-cond}
	\end{equation}
	
	Then $u(t)$ satisfies the following regularity properties.
	
	\begin{enumerate}
		\renewcommand{\labelenumi}{(\arabic{enumi})}
		
		\item \emph{(Regularity in the phase space)} It turns out that
		\begin{equation}
			(u,u')\in 
			C^{0}\left([0,+\infty),
			D(A^{\alpha_{0}})\times D(A^{\alpha_{1}})\right).
			\label{th:h:phsp}
		\end{equation}
		
		\item \emph{(Regularity of higher order derivatives up to
		$t=0$)} Let $m\geq 1$ be an integer such that $\alpha_{1}\geq
		(m-1)\gamma$.  Then the $m$-th time-derivative $u^{(m)}(t)$
		satisfies
		\begin{equation}
			u^{(m)}\in C^{0}\left([0,+\infty),
			D(A^{\alpha_{1}-(m-1)\gamma})\right).
			\label{th:h:0-m}
		\end{equation}

		\item \emph{(Regularity for $t>0$ when $0<\sigma<1$)} In 
		this regime it turns out that
		\begin{equation}
			u\in C^{\infty}((0,+\infty),D(A^{\alpha}))
			\quad\quad
			\forall\alpha\geq 0.
			\label{th:h:>0}
		\end{equation}
		
		\item \emph{(Regularity for $t>0$ when $\sigma\geq 1$)} In 
		this regime it turns out that
		\begin{equation}
			u^{(m)}\in C^{0}\left((0,+\infty),
			D(A^{\alpha_{0}+m(\sigma-1)})\right)
			\quad\quad
			\forall m\in\n.
			\label{th:h:>0'}
		\end{equation}
	\end{enumerate}
\end{thm}
\begin{rmk}\label{rmk:cont-dep}
	\begin{em}
		A careful inspection of the proofs reveals that the norm of
		the solution in the spaces appearing in~(\ref{th:h:phsp})
		through~(\ref{th:h:>0'}) depends continuously on the norm of
		initial data in $D(A^{\alpha_{0}})\times D(A^{\alpha_{1}})$.
		For example, in the case of (\ref{th:h:phsp}) this means that
		\begin{equation}
			\|u(t)\|_{D(A^{\alpha_{0}})}+\|u'(t)\|_{D(A^{\alpha_{1}})}
			\leq C_{1}\left(\|u_{0}\|_{D(A^{\alpha_{0}})}+
			\|u_{1}\|_{D(A^{\alpha_{1}})}\right) 
			\quad\quad
			\forall t\geq 0
			\label{est:semigroup}
		\end{equation}
		for a suitable constant
		$C_{1}=C_{1}(\delta,\sigma,\alpha_{0},\alpha_{1})$, while in
		in the case of (\ref{th:h:>0}) this means that for every
		$m\in\n$ one has that
		$$\left\|u^{(m)}(t)\right\|_{D(A^{\alpha})} \leq
		\frac{C_{2}}{t^{C_{3}}}\left(\|u_{0}\|_{D(A^{\alpha_{0}})}+
		\|u_{1}\|_{D(A^{\alpha_{1}})}\right) 
		\quad\quad
		\forall t> 0$$
		for suitable constants $C_{2}$ and $C_{3}$, both depending on 
		$\delta$, $\sigma$, $\alpha_{0}$, $\alpha_{1}$, $\alpha$ and 
		$m$. We spare the reader from these standard details.
	\end{em}
\end{rmk}

\begin{rmk}\label{rmk:homog}
	\begin{em}
		Statement~(1) of Theorem~\ref{thm:homog} and estimate
		(\ref{est:semigroup}) are equivalent to saying that equation
		(\ref{pbm:h-eqn}) generates a continuous semigroup in
		$D(A^{\alpha_{0}})\times D(A^{\alpha_{1}})$ provided that
		inequality (\ref{gap-cond}) is satisfied.  In addition, a
		simple inspection of the proof reveals that, when $A$ is
		unbounded, inequality (\ref{gap-cond}) is also a necessary
		condition for equation (\ref{pbm:h-eqn}) to generate a
		continuous semigroup in $D(A^{\alpha_{0}})\times
		D(A^{\alpha_{1}})$.
		
		Thus (\ref{gap-cond}) describes all possible values of the
		phase space gap $\alpha_{0}-\alpha_{1}$.  These values are
		represented as a function of $\sigma$ by the shaded region in
		the following picture.
		
		\begin{center}
			\definecolor{verdino}{RGB}{91,247,182}
			\definecolor{verdone}{RGB}{23,76,50}
			\psset{unit=5ex}
			\pspicture(-1,-2.8)(5,4.5)
			\pspolygon*[linecolor=verdino](1,1)(4,4)(4,-2)
			\psline[linewidth=.7\pslinewidth]{->}(-0.5,0)(4.5,0)
			\psline[linewidth=.7\pslinewidth]{->}(0,-2)(0,4.2)
			\psline[linecolor=verdone,linewidth=2\pslinewidth](0,1)(1,1)(4,4)
			\psline[linecolor=verdone,linewidth=2\pslinewidth](1,1)(4,-2)
			\psline[linestyle=dashed,linecolor=blue,linewidth=.7\pslinewidth](1,-0.1)(1,1)
			\psdots[linecolor=verdone,linewidth=2\pslinewidth](0,1)
			\rput(1,-0.4){$\frac{1}{2}$}
			\rput(-0.3,1){$\frac{1}{2}$}
			\rput(4.5,-0.4){$\sigma$}
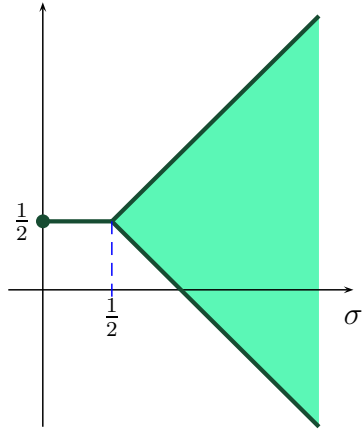
			\rput(1.75,-2.4){Admissible phase space gaps $\alpha_{0}-\alpha_{1}$ as a function of 
			$\sigma$}
			\endpspicture
		\end{center}
		
		We stress that for $0\leq\sigma\leq 1/2$ the only admissible
		value is $1/2$, which is typical of hyperbolic problems.  When
		$\sigma>1/2$ there is an interval of possible phase space
		gaps, centered in 1/2, which contains also negative values
		when $\sigma>1$.
		
		This implies that (\ref{pbm:eqn}) always generates a semigroup
		on $D(A^{1/2})\times H$, or more generally on
		$D(A^{\alpha+1/2})\times D(A^{\alpha})$, but for $\sigma>1/2$ there are
		always many other possible choices. Just to give some extremal 
		examples,
		\begin{itemize}
			\item  equation (\ref{pbm:h-eqn}) with $\sigma=1$ generates 
			a semigroup on $D(A)\times H$ or $H\times H$,
		
			\item equation (\ref{pbm:h-eqn}) with $\sigma=2$ generates a
			semigroup on $D(A^{2})\times H$ or $H\times D(A)$ (note
			that the latter has a negative phase space gap, namely the
			time-derivative is more regular than the function itself).
		\end{itemize}
	\end{em}
\end{rmk}

Now we proceed to the non-homogeneous case. The first question we 
address is the regularity of solutions. It is well-known that the 
space regularity of solutions to non-homogeneous linear equations 
depends both on the space and on the time regularity of the forcing 
term $f(t)$, in such a way that a higher time-regularity compensates 
a lower space-regularity. A typical example of this philosophy is 
Proposition~4.1.6 in~\cite{C-H}.

A full understanding of this interplay between time and space 
regularity in the case of equations with strong dissipation is 
probably an interesting problem, which could deserve further 
investigation. Here, for the sake of brevity, we limit ourselves to 
forcing terms with minimal space-regularity (just in $H$), and 
bounded with respect to time.

\begin{thm}[Non-homogeneous equation -- Regularity]\label{thm:nh-local}
	Let $H$ be a separable Hilbert
	space, and let $A$ be a self-adjoint nonnegative operator on $H$
	with dense domain $D(A)$.  Let $T>0$, and let $f\in
	L^{\infty}((0,T),H)$ be a bounded forcing term.
	
	For every $\sigma\geq 0$ and $\delta>0$, we consider the unique
	solution $u(t)$ of the non-homogeneous linear equation
	\begin{equation}
		u''(t)+2\delta A^{\sigma}u'(t)+Au(t)=f(t)
		\label{pbm:nh-eqn}
	\end{equation}
	in $[0,T]$, with null initial data
	\begin{equation}
		u(0)=0,
		\quad\quad
		u'(0)=0.
		\label{pbm:nh-data}
	\end{equation}
	
	Then $u(t)$ satisfies the following properties.
	\begin{enumerate}
		\renewcommand{\labelenumi}{(\arabic{enumi})}
		
		\item \emph{(Case $\sigma=0$)} In this regime it turns out 
		that
		$$u\in C^{0}([0,T],D(A^{1/2}))\cap C^{1}([0,T],H).$$
		
		\item \emph{(Case $0<\sigma<1$)} In this regime it turns out 
		that
		$$u\in C^{0}([0,T],D(A^{\min\{\sigma+1/2,1\}-\ep}))\cap 
		C^{1}([0,T],D(A^{\sigma-\ep}))
		\quad\quad
		\forall\ep\in(0,\sigma].$$
		
		\item \emph{(Case $\sigma\geq 1$)} In this regime it turns out 
		that
		\begin{equation}
			u\in C^{0}([0,T],D(A^{\sigma}))\cap 
			C^{1}([0,T],D(A^{\sigma-\ep}))
			\quad\quad
			\forall\ep\in(0,\sigma].
			\label{th:nhl:>1}
		\end{equation}
	\end{enumerate}
\end{thm}

\begin{rmk}\label{rmk:nh}
	\begin{em}
		The following two pictures sum up the conclusions of
		Theorem~\ref{thm:nh-local}.  The shaded regions represent the
		pairs $(\sigma,\alpha)$ for which it happens that $u$ is
		continuous with values in $D(A^{\alpha})$ and $u'$ is continuous with
		values in $D(A^{\alpha})$.  The boundary is included when
		represented by a continuous line or a dot, and excluded when
		represented by a dashed line.
		\begin{center}
			\hfill
			\definecolor{verdino}{RGB}{91,247,182}
			\definecolor{verdone}{RGB}{23,76,50}
			\psset{unit=5ex}
			
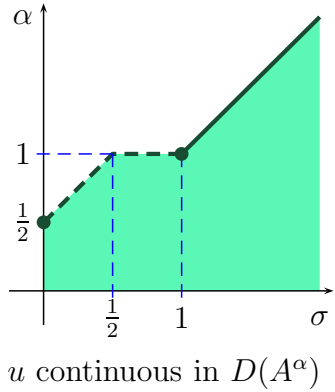
\pspicture(-1,-1.4)(4.5,4.5)
			\pspolygon*[linecolor=verdino](0,1)(1,2)(2,2)(4,4)(4,0)(0,0)
			\psline[linewidth=.7\pslinewidth]{->}(-0.5,0)(4.2,0)
			\psline[linewidth=.7\pslinewidth]{->}(0,-0.5)(0,4.2)
			\psline[linestyle=dashed,linecolor=verdone,linewidth=2\pslinewidth](0,1)(1,2)(2,2)
			\psline[linecolor=verdone,linewidth=2\pslinewidth](2,2)(4,4)
			\psline[linestyle=dashed,linecolor=blue,linewidth=.7\pslinewidth](2,-0.1)(2,2)
			\psline[linestyle=dashed,linecolor=blue,linewidth=.7\pslinewidth](1,-0.1)(1,2)
			\psline[linestyle=dashed,linecolor=blue,linewidth=.7\pslinewidth](-0.1,2)(1,2)
			\psdots[linecolor=verdone,linewidth=2\pslinewidth](0,1)(2,2)
			\rput(1,-0.4){$\frac{1}{2}$}
			\rput(2,-0.4){$1$}
			\rput(-0.3,1){$\frac{1}{2}$}
			\rput(-0.3,2){$1$}
			\rput(4,-0.4){$\sigma$}
			\rput(-0.3,4){$\alpha$}
			\rput(1.75,-1.2){$u$ continuous in $D(A^{\alpha})$}
			\endpspicture
			\hfill\hfill
			
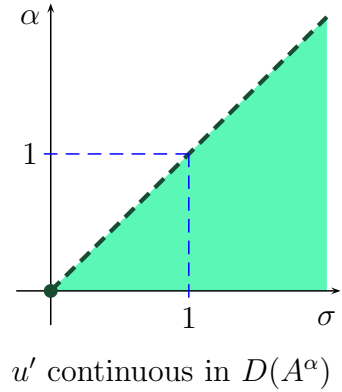
\pspicture(-1,-1.4)(4.5,4.5)
			\pspolygon*[linecolor=verdino](0,0)(4,4)(4,0)
			\psline[linewidth=.7\pslinewidth]{->}(-0.5,0)(4.2,0)
			\psline[linewidth=.7\pslinewidth]{->}(0,-0.5)(0,4.2)
			\psline[linestyle=dashed,linecolor=verdone,linewidth=2\pslinewidth](0,0)(4,4)
			\psline[linestyle=dashed,linecolor=blue,linewidth=.7\pslinewidth](2,-0.1)(2,2)
			\psline[linestyle=dashed,linecolor=blue,linewidth=.7\pslinewidth](-0.1,2)(2,2)
			\psdots[linecolor=verdone,linewidth=2\pslinewidth](0,0)
			\rput(2,-0.4){$1$}
			\rput(-0.3,2){$1$}
			\rput(4,-0.4){$\sigma$}
			\rput(-0.3,4){$\alpha$}
			\rput(1.75,-1.2){$u'$ continuous in $D(A^{\alpha})$}
			\endpspicture
			\hfill\mbox{}
		\end{center}
		
		As expected, $u(t)$ is always more regular than $u'(t)$.  One
		could call ``smoothing gap'' the difference between the
		exponents of the spaces where $u(t)$ and $u'(t)$ lie.  Thus
		the smoothing gap is equal to 1/2 (as usual in hyperbolic
		problems) for $\sigma\in[0,1/2]$, then it is equal to
		$1-\sigma$ for $\sigma\in(1/2,1)$, and finally it is 0 for
		$\sigma\geq 1$, when in any case $u(t)$ lies in the limit
		space $D(A^{\sigma})$ while $u'(t)$ does not.
	\end{em}
\end{rmk}

In the following result we investigate the dashed lines of the 
pictures above, showing that in those limit cases we have at least 
that $u(t)\in D(A^{\alpha})$ or $u'(t)\in D(A^{\alpha})$ for almost 
every time.

\begin{thm}[Non-homogeneous equation -- Limit cases]\label{thm:nh-limit}
	Let $H$ be a separable Hilbert space, and let $A$ be a
	self-adjoint nonnegative operator on $H$ with dense domain $D(A)$,
	let $T>0$, and let $f\in L^{2}((0,T),H)$.
	
	For every $\sigma\geq 0$ and $\delta>0$, we consider the unique
	solution $u(t)$ of the non-homogeneous linear 
	equation~(\ref{pbm:nh-eqn}) in $[0,T]$, with null initial 
	data~(\ref{pbm:nh-data}).
	
	Then $u(t)$ satisfies the following regularity properties, 
	depending on $\sigma$.
	\begin{enumerate}
		\renewcommand{\labelenumi}{(\arabic{enumi})}
		
		\item \emph{(Case $\sigma\geq 0$)} For every admissible value 
		of $\sigma$ it turns out that
		\begin{equation}
			u'\in L^{2}((0,T),D(A^{\sigma})).
			\label{th:nhl:u'2}
		\end{equation}
		
		\item \emph{(Case $\sigma\in[0,1]$)} In this regime it turns
		out that
		\begin{equation}
			u\in L^{2}\left((0,T),D(A^{\min\{\sigma+1/2,1\}})\right).
			\label{th:nhl:u2}
		\end{equation}
		
	\end{enumerate}
\end{thm}

The next result concerns the boundedness of solutions. We assume that the 
forcing term $f(t)$ is defined for every $t\geq 0$ and globally 
bounded, and we characterize the spaces where the solution is 
globally bounded. In this case we need to assume that the operator 
$A$ is coercive, because if not there are trivial counterexamples 
(just think to the case where $A$ is the null operator and $f$ is 
constant).

\begin{thm}[Non-homogeneous equation -- Global boundedness]\label{thm:nh-bound} 
	Let $H$ be a separable Hilbert space, and let $A$ be a
	self-adjoint nonnegative operator on $H$ with dense domain $D(A)$.
	Let us assume that $A$ is coercive, namely there exists a constant
	$\nu>0$ such that $\langle Au,u\rangle\geq\nu|u|^{2}$ for every
	$u\in D(A)$.  Let $f\in L^{\infty}((0,+\infty),H)$ be a globally
	bounded forcing term.
	
	For every $\sigma\geq 0$ and $\delta>0$, we consider the unique
	solution $u(t)$ of the non-homogeneous linear
	equation~(\ref{pbm:nh-eqn}) in $[0,+\infty)$, with null initial
	data~(\ref{pbm:nh-data}).
	
	Then $u(t)$ is bounded in the following spaces, depending on 
	$\sigma$.
	\begin{enumerate}
		\renewcommand{\labelenumi}{(\arabic{enumi})}
		
		\item \emph{(Case $\sigma=0$)} In this regime it turns out that
		$$(u(t),u'(t))
		\mbox{ is bounded in }
		D(A^{1/2})\times H.$$
		
		\item \emph{(Case $0<\sigma<1$)} In this regime it turns out that
		$$
		(u(t),u'(t))
		\mbox{ is bounded in }
		D(A^{\min\{\sigma+1/2,1\}-\ep})\times D(A^{\sigma-\ep})
		\quad
		\forall\ep\in(0,\sigma].
		$$
		
		\item \emph{(Case $\sigma= 1$)} In this regime it turns out that
		$$
		(u(t),u'(t))
		\mbox{ is bounded in }
		D(A)\times D(A^{1-\ep})
		\quad
		\forall\ep\in(0,1].
		$$
		
		\item \emph{(Case $\sigma> 1$)} In this regime it turns out that
		$$
		(u(t),u'(t))
		\mbox{ is bounded in }
		D(A^{1-\ep})\times D(A^{\sigma-\ep})
		\quad
		\forall\ep\in(0,1].
		$$
		
	\end{enumerate}
\end{thm}

\begin{rmk}
	\begin{em}
		In the proof of Theorems~\ref{thm:nh-local},
		\ref{thm:nh-limit} and \ref{thm:nh-bound} we actually show
		also that the norm of the solution (in the given spaces)
		depends continuously on the norm of $f$.  Just to give an
		example, in the case of (\ref{th:nhl:>1}) this means that
		$$\|u(t)\|_{D(A^{\sigma})}\leq
		C_{1}\|f\|_{L^{\infty}((0,T),H)} \quad\quad
		\forall t\in [0,T],$$
		for a suitable constant $C_{1}=C_{1}(\delta,\sigma)$, and
		$$\|u'(t)\|_{D(A^{\sigma-\ep})}\leq C_{2}\|f\|_{L^{\infty}((0,T),H)}
		\quad\quad
		\forall t\in [0,T],$$
		for a suitable constant $C_{2}=C_{2}(\delta,\sigma,\ep)$.
	\end{em}
\end{rmk}

\begin{rmk}\label{rmk:bounded}
	\begin{em}
		The following two pictures sum up the conclusions of
		Theorem~\ref{thm:nh-bound}, namely the spaces where $u(t)$ and
		$u'(t)$ are globally bounded when the forcing term is globally
		bounded.
		
		\begin{center}
			\hfill
			\definecolor{verdino}{RGB}{91,247,182}
			\definecolor{verdone}{RGB}{23,76,50}
			\psset{unit=5ex}
			
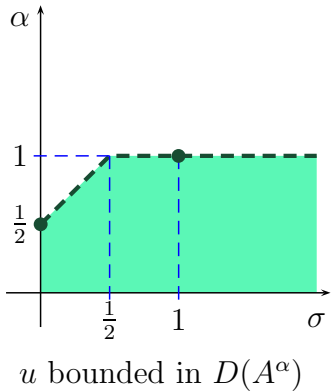
\pspicture(-1,-1.4)(4.5,4.5)
			\pspolygon*[linecolor=verdino](0,1)(1,2)(2,2)(4,2)(4,0)(0,0)
			\psline[linewidth=.7\pslinewidth]{->}(-0.5,0)(4.2,0)
			\psline[linewidth=.7\pslinewidth]{->}(0,-0.5)(0,4.2)
			\psline[linestyle=dashed,linecolor=verdone,linewidth=2\pslinewidth](0,1)(1,2)(4,2)
			\psline[linestyle=dashed,linecolor=blue,linewidth=.7\pslinewidth](2,-0.1)(2,2)
			\psline[linestyle=dashed,linecolor=blue,linewidth=.7\pslinewidth](1,-0.1)(1,2)
			\psline[linestyle=dashed,linecolor=blue,linewidth=.7\pslinewidth](-0.1,2)(1,2)
			\psdots[linecolor=verdone,linewidth=2\pslinewidth](0,1)(2,2)
			\rput(1,-0.4){$\frac{1}{2}$}
			\rput(2,-0.4){$1$}
			\rput(-0.3,1){$\frac{1}{2}$}
			\rput(-0.3,2){$1$}
			\rput(4,-0.4){$\sigma$}
			\rput(-0.3,4){$\alpha$}
			\rput(1.75,-1.2){$u$ bounded in $D(A^{\alpha})$}
			\endpspicture
			\hfill\hfill
			
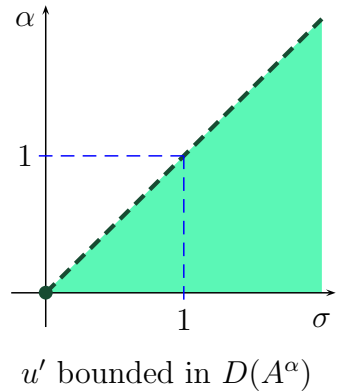
\pspicture(-1,-1.4)(4.5,4.5)
			\pspolygon*[linecolor=verdino](0,0)(4,4)(4,0)
			\psline[linewidth=.7\pslinewidth]{->}(-0.5,0)(4.2,0)
			\psline[linewidth=.7\pslinewidth]{->}(0,-0.5)(0,4.2)
			\psline[linestyle=dashed,linecolor=verdone,linewidth=2\pslinewidth](0,0)(4,4)
			\psline[linestyle=dashed,linecolor=blue,linewidth=.7\pslinewidth](2,-0.1)(2,2)
			\psline[linestyle=dashed,linecolor=blue,linewidth=.7\pslinewidth](-0.1,2)(2,2)
			\psdots[linecolor=verdone,linewidth=2\pslinewidth](0,0)
			\rput(2,-0.4){$1$}
			\rput(-0.3,2){$1$}
			\rput(4,-0.4){$\sigma$}
			\rput(-0.3,4){$\alpha$}
			\rput(1.75,-1.2){$u'$ bounded in $D(A^{\alpha})$}
			\endpspicture
			\hfill\mbox{}
		\end{center}
		
		Comparing with Remark~\ref{rmk:nh}, we see that the regularity
		and boundedness diagrams of $u'(t)$ coincide for all
		$\sigma\geq 0$, while the regularity and boundedness diagrams
		of $u(t)$ coincide only for $\sigma\in[0,1]$.  In other words,
		for $\sigma>1$ the solution is more regular, but the estimates
		in the stronger norms diverge as $t\to+\infty$.  When
		$\sigma>1$ and $1\leq\alpha\leq\sigma$ one can obtain
		estimates of the form
		\begin{equation}
			|A^{\alpha}u(t)|\leq
			C_{\alpha,\sigma}\|f\|_{L^{\infty}((0,+\infty),H)}\cdot
			\left\{
			\begin{array}{ll}
				t^{(\alpha-1)/(\sigma-1)} & 
				\mbox{if }1<\alpha\leq\sigma, \\
				\noalign{\vspace{0.5ex}}
				\log(1+t) & \mbox{if }\alpha=1.
			\end{array}
			\right.
			\label{th:growth}
		\end{equation}
		
		We refer to Remark~\ref{rmk:lemma-growth} for further details.
		We observe also that $\sigma=1$ plays a special r\^{o}le in
		the boundedness diagrams, being the unique exponent for which
		there is global boundedness of $u(t)$ in $D(A)$.
		
	\end{em}
\end{rmk}

\begin{rmk}
	\begin{em}
		Regularity and boundedness properties of $u''(t)$, or more
		generally of further time-derivatives of $u(t)$, can be easily
		deduced from the regularity and boundedness properties of the
		other three terms in equation~(\ref{pbm:eqn}).  Thus we see
		that there is no value of $\sigma$ for which a forcing term
		$f\in L^{\infty}((0,T),H)$ is enough to guarantee that all
		terms in the left-hand side of (\ref{pbm:eqn}) make sense
		individually as elements of $H$.  Therefore, solutions are
		always to be intended as weak solutions.  We observe also that
		the terms $Au(t)$ and $A^{\sigma}u'(t)$ are in the same spaces
		if and only if $1/2\leq\sigma<1$.
		
	\end{em}
\end{rmk}

For the next result we consider equation~(\ref{pbm:eqn}) with a
forcing term $f\in L^{\infty}(\re,H)$.  We look for a solution $u(t)$
which is bounded in some sense for every $t\in\re$.  As before, we
restrict ourselves to coercive operators, because if not the existence
of such a solution is in general false (a simple example being when
$A=0$ and $f$ is constant).  Whenever $0\leq \sigma\leq 1$, the
homogeneous system is exponentially damped in the standard energy
space, and in this case the classical result is existence and
uniqueness of a bounded solution on the line which attracts
exponentially all solutions as $t$ tends to $+\infty$.  On the other
hand when $\sigma> 1$ and $A$ is for instance diagonal and unbounded,
exponential damping is no longer satisfied.  However in this case a slight
modification of the tools used to prove Theorem~\ref {thm:nh-bound}
will give us not only the existence of a bounded solution on the line, but 
the properties on the line that we had previously on the half-line.

\bigskip  

In order to state properly the next result on the line, we need to
introduce briefly some notation.  First, given any Hilbert space $X$
and any possibly unbounded closed interval $J$, the space of all
functions $f:J\to X$ which are continuous and bounded, endowed with
the uniform norm on $J$, will be denoted from now on by
$C^{0}_{b}(J,X)$.  Then, following Bochner's definition, a function
$f\in C^{0}_b(\re, X)$ will be called almost periodic with values in $X$
iff the set of translates 
$$\bigcup_{\alpha\in \re} \{ f(t+\alpha)\}$$ 
is precompact in the space $C^{0}_b(\re, X)$.  The Banach space of
such functions, endowed by the topology of $ C^{0}_b(\re, X)$, is
denoted by $AP(\re, X)$.  Any $f\in AP(\re, X)$ can be represented by
a formal expansion on the complexified extension of $X$ of the form
$$f \sim \sum_ {j\in \n} f_ j e^{\mu_j t},$$ 
where the real numbers $\mu_j$ are defined by the property that they
belong to those numbers $\mu$ for which 
$$ \lim_{T\rightarrow\infty}\frac{1}{T}\int_0^Tf (t) e^{-\mu t} dt \not = 0.$$

The set of all such real numbers $\mu$ depends on $f$, is always
countable and is denoted by $\exp(f)$.  A very important property of
almost periodic function is the following: if $X$, $Y$ are two real
Hilbert spaces and ${\cal C} \in L (C_b(\re, X), C_b(\re, Y))$ is a
bounded linear operator then it is immediate, using Bochner's
definition, to see that 
$$\forall f \in AP(\re, X), \quad {\cal C} f
\in AP(\re, Y) \quad \hbox {with}\quad \exp({\cal C}
f)\subseteq \exp(f).$$

For more details on these questions, the construction of the
mean-value, the proof that the set $\exp(f)$ is countable and the
exact meaning of the formal expansion, we refer to \cite{LZ}.  A
typical almost periodic numerical function is the sum of two periodic
functions with incommensurable periods.  Such objects often appear in
the mechanics of vibrating systems, and sometimes infinite sums
naturally impose their presence, for instance when studying the energy
conservative vibrations of continuous media.  For classical
applications and historical comments, cf.  e.g. \cite{AP}.

\begin{thm}[Existence and properties of the bounded solution]\label{thm:glob-bound}
	Let $H$ and $A$ be as in Theorem~\ref{thm:nh-bound} (in particular
	the operator $A$ is assumed to be coercive).  Let $f\in
	L^{\infty}(\re,H)$ be a globally bounded forcing term.
	
	Then, for every $\sigma\geq 0$ and every $\delta>0$, equation
	(\ref{pbm:eqn}) admits a unique global solution which is
	continuous and bounded in the energy space $D(A^{1/2})\times H$.
	This solution is strongly asymptotic in the energy space to any
	solution with initial data in $D(A^{1/2})\times H$.
	
	Moreover, this solution satisfies the following regularity and
	boundedness properties.
	\begin{enumerate}
		\renewcommand{\labelenumi}{(\arabic{enumi})} 
		
		\item  The pair $(u(t),u'(t))$ is continuous and bounded in 
		the same spaces as those of Theorem~\ref{thm:nh-bound}, namely
		\begin{itemize}
			\item  in $D(A^{1/2})\times H$ if $\sigma=0$,
		
			\item in $D(A^{\min\{\sigma+1/2,1\}-\ep})\times
			D(A^{\sigma-\ep})$ for every $\ep\in(0,\sigma]$ if
			$0<\sigma<1$,
			
			\item in $D(A)\times D(A^{1-\ep})$ for every $\ep\in(0,1]$
			if $\sigma=1$,
		
			\item in $D(A^{1-\ep})\times D(A^{\sigma-\ep})$ for every
			$\ep\in(0,1]$ if $\sigma>1$.
			
		\end{itemize}
		
		\item If in addition $f$ is almost periodic with values in
		$H$, then $(u, u')$ is almost periodic with values in the
		spaces mentioned above with
		$\exp(u)\subseteq \exp(f)$.  Finally, if $f$ is
		periodic the bounded solution is periodic as well, with the
		same minimal period.  If $\sigma>1$, the periodic solution is
		continuous and bounded also in the limit space $D(A)$.
	
	\end{enumerate}
	
\end{thm}

\medskip  

Finally, we show that all previous results are optimal.  Note that in
the counterexamples below we always produce forcing terms which are
not just bounded, but also continuous.  This shows that a
time-continuous external force does not make the solution more
space-regular than a time-bounded external force.

\begin{thm}[Counterexamples]\label{thm:optimal} 
	Let $H$ be a Hilbert space, and let $A$ be a linear operator on
	$H$ with domain $D(A)$.  Let us assume that the spectrum of $A$
	contains an unbounded sequence of positive eigenvalues. 
	Then we have the following conclusions.
	\begin{enumerate}
		\renewcommand{\labelenumi}{(\arabic{enumi})} 
		
		\item \emph{(Case $\sigma=0$)} For every sequence
		$\{t_{n}\}\subseteq(0,+\infty)$ there exists $f\in
		C^{0}_{b}([0,+\infty),H)$ such that the unique global solution $u(t)$
		of problem (\ref{pbm:nh-eqn})--(\ref{pbm:nh-data}) satisfies
		\begin{equation}
			u(t_{n})\not\in D(A^{1/2+\ep})
			\quad\quad
			\forall\ep>0,\quad\forall n\in\n,
			\label{th:1a}
		\end{equation}
		\begin{equation}
			u'(t_{n})\not\in D(A^{\ep})
			\quad\quad
			\forall\ep>0,\quad\forall n\in\n.
			\label{th:1b}
		\end{equation}
	
		\item \emph{(Case $0<\sigma<1$)} For every $\sigma$ in this
		range, and every sequence $\{t_{n}\}\subseteq(0,+\infty)$,
		there exists $f\in C^{0}_{b}([0,+\infty),H)$ such that the
		unique global solution $u(t)$ of problem
		(\ref{pbm:nh-eqn})--(\ref{pbm:nh-data}) satisfies
		$$
		u(t_{n})\not\in D(A^{\min\{\sigma+1/2,1\}})
		\quad\quad\quad
		\forall n\in\n,
		$$
		$$
		u'(t_{n})\not\in D(A^{\sigma})
		\quad\quad\quad
		\forall n\in\n.
		$$
	
		\item \emph{(Case $\sigma\geq 1$)} For every $\sigma$ in this
		range, and every sequence $\{t_{n}\}\subseteq(0,+\infty)$,
		there exists $f\in C^{0}_{b}([0,+\infty),H)$ such that the
		unique global solution $u(t)$ of problem
		(\ref{pbm:nh-eqn})--(\ref{pbm:nh-data}) satisfies
		\begin{equation}
			u(t)\not\in D(A^{\sigma+\ep})
			\quad\quad
			\forall\ep>0,\quad\forall t>0,
			\label{th:4a}
		\end{equation}
		\begin{equation}
			u'(t_{n})\not\in D(A^{\sigma})
			\quad\quad\quad
			\forall n\in\n.
			\label{th:4b}
		\end{equation}
	
		\item \emph{(Case $\sigma> 1$ -- Unboundedness)} For every
		$\sigma$ in this range, there exists $f\in
		C^{0}_{b}([0,+\infty),H)$, and a sequence $t_{n}\to +\infty$, such
		that the unique global solution $u(t)$ of problem
		(\ref{pbm:nh-eqn})--(\ref{pbm:nh-data}) satisfies
		\begin{equation}
			\lim_{n\to +\infty}|Au(t_{n})|= +\infty.
			\label{th:5}
		\end{equation}
	\end{enumerate}
\end{thm}

\begin{rmk}
	\begin{em}
		In most counterexamples we obtained that $u$ or $u'$ does not
		belong to a given space for a sequence of times.  This cannot
		be improved by showing the same for all times, at least when
		the counterexample refers to the dashed lines in the
		continuity diagrams of Remark~\ref{rmk:nh}.  Indeed, in those
		cases we already observed that $u$ or $u'$ lie in the spaces
		of the dashed lines for almost every time (see
		Theorem~\ref{thm:nh-limit}).
		
		Finally, we point out that a careful inspection of the proof
		reveals that in our counterexample the norm in~(\ref{th:5})
		blows-up logarithmically with respect to $t_{n}$.  This proves
		the optimality of our growth estimates~(\ref{th:growth}) in
		the limit case $\alpha=1$.
		
	\end{em}
\end{rmk}

\begin{rmk}\label{rmk:final}
	\begin{em} 
		
		The regularity of the bounded solution of
		Theorem~\ref{thm:glob-bound} is also optimal.  Indeed it is
		enough to consider a forcing term which is 0 for $t\leq 0$ and
		equal to one of the forcing terms of Theorem~\ref{thm:optimal}
		for $t\geq 0$.  It is clear that with this choice the bounded
		solution is~0 for $t\leq 0$ (in particular $u(0)=u'(0)=0$),
		and equal to the corresponding solution in
		Theorem~\ref{thm:optimal} for positive times.  This yields the
		required regularity loss.
		
		It is possible to show that our regularity statements are
		optimal also in the case of periodic forcing terms.  The idea
		is to take the above forcing terms, which are always defined in
		a bounded time interval, and to extend them by periodicity.
		For the sake of shortness, we do not work out the details.
	\end{em}
\end{rmk}

\setcounter{equation}{0}

\section{Proofs in the homogeneous case}\label{sec:proofs-autonomous}

\subsection{Notation}\label{sec:notation}

Throughout this Section, $H$ denotes a separable Hilbert space, and 
$A$ is any self-adjoint nonnegative linear operator on $H$.  According
to the spectral theorem (see for example Theorem~VIII.4
of~\cite{reed}), there exist a measure space $(M,\mu)$ with $\mu$ a
finite measure, a unitary operator $\Psi:H\to L^{2}(M,\mu)$, and a
real valued nonnegative function $\lambda(\xi)$, defined for almost
every $\xi\in M$, such that
\begin{itemize}
	\item  $u\in D(A)$ if and only if $\lambda(\xi)[\Psi(u)](\xi)\in 
	L^{2}(M,\mu)$,

	\item $[\Psi(Au)](\xi)=\lambda(\xi)[\Psi(u)](\xi)$ for almost
	every $\xi\in M$.
\end{itemize}

In other words, $H$ can be identified with $L^{2}(M,\mu)$, and under
this identification the operator $A$ becomes the multiplication
operator by $\lambda(\xi)$.  In the sequel, for the sake of
simplicity, we write $\ut(\xi)$ instead of $[\Psi(u)](\xi)$.  We also
think of $\ut(\xi)$ as the ``component'' of $u$ with respect to some
$\xi\in M$.  One can work with these ``components'' exactly as with
Fourier series or Fourier transforms.  For example, for every
$\alpha\geq 0$ it is true that
$$D(A^{\alpha}):=\left\{u\in H:
[\laxi]^{\alpha}\cdot\ut(\xi)\in L^{2}(M,\mu)\right\},$$
and the components of $A^{\alpha}u$ are just
$[\laxi]^{\alpha}\cdot\ut(\xi)$.  Moreover, $u(t)$ turns out to be a
solution to (\ref{pbm:eqn}) if and only if the components $\ut(t,\xi)$
of $u(t)$ and the components $\ft(t,\xi)$ of $f(t)$ satisfy the
ordinary differential equation
\begin{equation}
	\ut\,''(\xi,t)+2\delta\,[\laxi]^{\sigma}\cdot\ut\,'(t,\xi)+\lambda(\xi) 
	\ut(t,\xi)=\ft(t,\xi)
	\label{eqn:uxi}
\end{equation}
for almost every $\xi\in M$, where as usual primes denote
differentiation with respect to $t$. 

In the proof of our regularity results we can always assume that the
operator $A$ is unbounded, because if not all regularity
statements in Theorems~\ref{thm:homog}, \ref{thm:nh-local},
and~\ref{thm:nh-bound} are trivial.  We can also assume that
multipliers $\laxi$ are large enough.  Indeed
the linearity of the equation implies that, given any threshold
$\Lambda>0$, every solution of (\ref{pbm:eqn})--(\ref{pbm:data}) is
the sum of two solutions, one corresponding to components with respect
multipliers $\laxi<\Lambda$ (and this solution
is regular), and one corresponding to components with respect to
multipliers $\laxi\geq\Lambda$.\medskip

In order to simplify the notation, in the sequel we always write ``for
every $\xi\in M$'' instead of ``for almost every $\xi\in M$'', and
when we write a supremum over $M$ we actually mean an essential
supremum. 

\subsection{Roots of the characteristic polynomial}\label{sec:roots}

The behavior of solutions of (\ref{eqn:uxi}) depends on the roots of
the characteristic polynomial
\begin{equation}
	x^{2}+2\delta\laxi^{\sigma}x+\laxi.
	\label{char-pol-xi}
\end{equation}

We call the roots $-x_{1}(\xi)$ and $-x_{2}(\xi)$ in order to
emphasize that they are negative real numbers, or complex numbers with
negative real part.  The asymptotic behavior of the roots as $\laxi\to
+\infty$ is different in the following three regimes.
\begin{itemize}
	\item \emph{Subcritical dissipation}.  If $0\leq\sigma<1/2$, or
	$\sigma=1/2$ and $\delta<1$, then for $\laxi$ large enough the
	roots of (\ref{char-pol-xi}) are complex conjugate numbers of the
	form
	\begin{equation}
		-x_{1}(\xi)=-a(\xi)+ib(\xi),
		\quad\quad\quad
		-x_{2}(\xi)=-a(\xi)-ib(\xi),
		\label{roots:sub-c}
	\end{equation}
	with
	$$a(\xi):=\delta\laxi^{\sigma}, \quad\quad\quad
	b(\xi):=\left(\laxi-\delta^{2}\laxi^{2\sigma}\right)^{1/2}.  $$
	
	As $\laxi\to +\infty$ it turns out that 
	$$b(\xi)\sim\left\{
	\begin{array}{ll}
		\laxi^{1/2}  &  \mbox{if $\sigma<1/2$,}  \\
		\noalign{\vspace{1ex}}
		\laxi^{1/2}(1-\delta^{2})^{1/2}  &  \mbox{if $\sigma=1/2$.}
	\end{array}
	\right.$$

	\item \emph{Critical dissipation}.  If $\sigma=1/2$ and
	$\delta=1$, then for every $\laxi\geq 0$ the characteristic
	polynomial (\ref{char-pol-xi}) has a unique root
	\begin{equation}
		-x_{1}(\xi)=-x_{2}(\xi)=-\laxi^{1/2} 
		\label{roots:c}
	\end{equation}
	with multiplicity 2.

	\item \emph{Supercritical dissipation}.  If $\sigma>1/2$, or
	$\sigma=1/2$ and $\delta>1$, then for $\laxi$ large enough the
	characteristic polynomial (\ref{char-pol-xi}) has two distinct
	real roots
	\begin{equation}
		\begin{array}{c}
			-x_{1}(\xi)= -\delta\laxi^{\sigma}-
			\left(\delta^{2}\laxi^{2\sigma}-\laxi\right)^{1/2},  \\
			\noalign{\vspace{1ex}}
			-x_{2}(\xi)= -\delta\laxi^{\sigma}+
			\left(\delta^{2}\laxi^{2\sigma}-\laxi\right)^{1/2},
		\end{array}
		\label{roots:super-c}
	\end{equation}
	so that for $\sigma>1/2$ it turns out that
	$$x_{1}(\xi)\sim 2\delta\laxi^{\sigma}
	\quad\quad\mbox{and}\quad\quad
	x_{2}(\xi)\sim \frac{1}{2\delta}\laxi^{1-\sigma}$$
	as $\laxi\to +\infty$, while for $\sigma=1/2$ it turns out that 
	$$x_{1}(\xi)=\left(\delta+\sqrt{\delta^{2}-1}\right)\laxi^{1/2}
	\quad\quad\mbox{and}\quad\quad
	x_{2}(\xi)=\left(\delta-\sqrt{\delta^{2}-1}\right)\laxi^{1/2}.$$
\end{itemize}

\subsection{A basic lemma}

Every solution of (\ref{eqn:uxi}) with $f\equiv 0$ can be written as a finite
sum of terms of the form $z_{0}(\xi)\cdot c(\xi)\cdot g(t,\xi)$, where
$z_{0}(\xi)$ is one of the initial conditions, $c(\xi)$ is a
coefficient depending only on the roots of the characteristic
polynomial~(\ref{char-pol-xi}), and $g(t,\xi)$ is one of the
fundamental solutions of the same equation (depending once again on
the roots of the characteristic polynomial).  For this reason we
investigate the regularity of these objects.

\begin{lemma}\label{lemma:main-h}
	Let $H$, $A$, $(M,\mu)$, $\laxi$ be as in
	section~\ref{sec:notation}.  Let $z_{0}(\xi)$ and $c(\xi)$ be two
	measurable real function on $M$, and let $g(t,\xi)$ be a function
	defined in $[0,+\infty)\times M$ which is of class $C^{\infty}$ with
	respect to $t$ and measurable with respect to $\xi$.
	
	Let us assume that there exists $P\geq 0$ such that 
	\begin{equation}
		S_{P}:=\int_{M}\laxi^{2P}z_{0}(\xi)^{2}\,d\mu(\xi)<+\infty,
		\label{hp:sp}
	\end{equation}
	and there exists $Q\geq 0$ such that 
	\begin{equation}
		C_{Q}:=\sup_{\xi\in M}\laxi^{Q}|c(\xi)|<+\infty.
		\label{hp:cq}
	\end{equation}	
	
	Then the following statements hold true.

	\begin{enumerate}
		\renewcommand{\labelenumi}{(\arabic{enumi})} 
		
		\item  Let us assume that there exists a constant $G_{0}$ 
		such that
		\begin{equation}
			|g(t,\xi)|\leq G_{0}
			\quad\quad
			\forall (t,\xi)\in [0,+\infty)\times M.
			\label{hp:G0}
		\end{equation}
		
		Then the product
		$$\zt(t,\xi):=z_{0}(\xi)\cdot c(\xi)\cdot g(t,\xi) 
		\quad\quad
		\forall (t,\xi)\in [0,+\infty)\times M$$
		defines a function $\zt$ which corresponds, under the usual
		identification of $L^{2}(M,\mu)$ with $H$, to a function $z\in
		C^{0}([0,+\infty),D(A^{P+Q}))$.
	
		\item More generally, let us assume that there exist an
		integer $m\geq 0$, and real constants $R$ (not necessarily
		positive) and $G_{m}\geq 0$ such that $P+Q-mR\geq 0$ and the
		$m$-th time-derivative $g^{(m)}(t,\xi)$ of $g(t,\xi)$
		satisfies 
		\begin{equation}
			|g^{(m)}(t,\xi)|\leq G_{m}\laxi^{mR} 
			\quad\quad
			\forall (t,\xi)\in [0,+\infty)\times M.
			\label{hp:Gm}
		\end{equation}
		
		Then $z$ is $m$ times differentiable (for example with values
		in $H$), and actually its $m$-th time-derivative $z^{(m)}(t)$
		satisfies
		$$z^{(m)}\in C^{0}\left([0,+\infty),D(A^{P+Q-mR})\right).$$
	
		\item Let us assume that there exist a measurable positive
		function $\eta(\xi)$ in $M$, a sequence
		$\{\Gamma_{m}\}\subseteq[0,+\infty)$, and constants $R\in\re$,
		$S>0$, and $M_{S}>0$ such that 
		\begin{equation}
			|g^{(m)}(t,\xi)|\leq
			\Gamma_{m}\laxi^{mR}e^{-\eta(\xi)t} 
			\quad\quad
			\forall (m,t,\xi)\in\n\times[0,+\infty)\times M,
			\label{hp:Gamma-m}
		\end{equation}
		and
		\begin{equation}
			\laxi^{S}\leq M_{S}\eta(\xi)
			\quad\quad
			\forall\xi\in M.
			\label{hp:ms}
		\end{equation}
		
		Then it turns out that $z\in
		C^{\infty}((0,+\infty),D(A^{\alpha}))$ for every $\alpha\geq 0$.
	\end{enumerate}
	
\end{lemma} 

\paragraph{\textmd{\textit{Proof}}}

\subparagraph{\textmd{\textit{Statement (1)}}}

Due to (\ref{hp:cq}) and (\ref{hp:G0}) it turns out that
\begin{eqnarray}
	\laxi^{2P+2Q}\cdot|\zt(t,\xi)|^{2} & = &
	\laxi^{2P}z_{0}(\xi)^{2}\cdot\laxi^{2Q}|c(\xi)|^{2}
	\cdot|g(t,\xi)|^{2}
	\nonumber  \\
	\noalign{\vspace{1ex}}
	 & \leq & C_{Q}^{2}\cdot G_{0}^{2}\cdot\laxi^{2P}z_{0}(\xi)^{2}.
	\label{est:pre-leb}
\end{eqnarray}

Therefore, from (\ref{hp:sp}) it follows that
$$\int_{M}\laxi^{2P+2Q}\cdot|\zt(t,\xi)|^2\,d\mu(\xi)\leq
C_{Q}^{2}\cdot G_{0}^{2}\cdot S_{P} <+\infty,$$
which is equivalent to saying that $z(t)\in D(A^{P+Q})$ for every
$t\geq 0$.  Since (\ref{est:pre-leb}) is uniform in time, the
continuity of $A^{P+Q}z(t)$ with respect to $t$ follows from the
continuity of $g(t,\xi)$ with respect to $t$ and Lebesgue's theorem.

\subparagraph{\textmd{\textit{Statement (2)}}}

The $m$-th time-derivative $\zt^{(m)}(t,\xi)$ of
$\zt(t,\xi)$ exists because $g(t,\xi)$ is $m$ times
differentiable.  If (\ref{hp:Gm}) holds true, then
\begin{eqnarray*}
	\laxi^{2P+2Q-2mR}\cdot|\zt^{(m)}(t,\xi)|^{2} & = &
	\laxi^{2P}z_{0}(\xi)^{2}\cdot
	\laxi^{2Q}c(\xi)^{2}\cdot\laxi^{-2mR}|g^{(m)}(t,\xi)|^{2} \\
	\noalign{\vspace{1ex}}
	 & \leq & C_{Q}^{2}\cdot
	G_{m}^{2}\cdot\laxi^{2P}z_{0}(\xi)^{2}.
\end{eqnarray*}

Therefore, from (\ref{hp:sp}) it follows that
$$\int_{M}\laxi^{2P+2Q-2mR}\cdot|\zt^{(m)}(t,\xi)|^{2}\,d\mu(\xi)\leq
C_{Q}^{2}\cdot G_{m}^{2}\cdot S_{P} <+\infty,$$
which is equivalent to saying that $z^{(m)}(t)\in D(A^{P+Q-mR})$.  The
continuity with values in the same space follows from the continuity
of $g^{(m)}(t,\xi)$ with respect to $t$ and Lebesgue's theorem as in
the first statement.

\subparagraph{\textmd{\textit{Statement (3)}}}

For every $\beta>0$ there exists a constant $K_{\beta}$ such that
$e^{-x}\leq K_{\beta}x^{-\beta}$ for every $x>0$. Therefore 
assumption (\ref{hp:Gamma-m}) implies that 
$$\laxi^{-2mR}\cdot\left|g^{(m)}(t,\xi)\right|^{2}\leq
\Gamma_{m}^{2}e^{-2\eta(\xi)t}\leq
\Gamma_{m}^{2}\cdot K_{\beta}^{2}\cdot
\frac{1}{[\eta(\xi)t]^{2\beta}}
\quad\quad
\forall(t,\xi)\in(0,+\infty)\times M.$$

Keeping (\ref{hp:ms}) into account, it follows that
$$\laxi^{-2mR+2\beta S}\cdot\left|g^{(m)}(t,\xi)\right|^{2}\leq
\Gamma_{m}^{2}\cdot K_{\beta}^{2}\cdot
\left(\frac{\laxi^{S}}{\eta(\xi)}\right)^{2\beta}
\frac{1}{t^{2\beta}}\leq
\frac{\Gamma_{m}^{2}\cdot K_{\beta}^{2}\cdot M_{S}^{2\beta}}{t^{2\beta}},$$
and finally
\begin{eqnarray*}
	\laxi^{2P+2Q-2mR+2\beta S}\cdot
	\left|\zt^{(m)}(t,\xi)\right|^{2} & = & 
	\laxi^{2Q}c(\xi)^{2}\cdot\laxi^{-2mR+2\beta S}
	\left|g^{(m)}(t,\xi)\right|^{2} \\
	   &  &  \mbox{}\cdot\laxi^{2P}z_{0}(\xi)^{2} \\
	   \noalign{\vspace{1ex}}
	 & \leq &  C_{Q}^{2}\cdot 
	 \frac{\Gamma_{m}^{2}\cdot K_{\beta}^{2}\cdot M_{S}^{2\beta}}{t^{2\beta}}
	 \cdot\laxi^{2P}z_{0}(\xi)^{2}
\end{eqnarray*}
for every $t>0$ and every $\xi\in M$.  As in the first two statements,
we integrate with respect to $\xi$ and what we obtain is equivalent to
saying that $z^{(m)}(t)\in D(A^{P+Q-mR+\beta S})$ for every $t>0$.
Also the continuity of $z^{(m)}(t)$ with values in the same space
follows as in the first two statements.

Since $S>0$, the exponent $P+Q-mR+\beta S$ can be made arbitrarily
large by choosing $\beta$ large enough, and this is enough to
establish that $z\in C^{\infty}((0,+\infty),D(A^{\alpha}))$ for every
$\alpha\geq 0$.\qed \medskip
\subsection{Proof of Theorem~\ref{thm:homog}}

Let $\ut_{0}(\xi)$ and $\ut_{1}(\xi)$ denote the components of $u_{0}$
and $u_{1}$.  It is easy to express the components $\ut(t,\xi)$ of the
solution in terms of $\ut_{0}(\xi)$, $\ut_{1}(\xi)$, and of the roots
of the characteristic polynomial (\ref{char-pol-xi}).  We distinguish
three cases depending on the discriminant of~(\ref{char-pol-xi}).

\subparagraph{\textmd{\emph{Supercritical dissipation}}}

Let us consider the case where $\sigma>1/2$, or $\sigma=1/2$ and
$\delta>1$, so that for $\lambda$ large enough the roots of the
characteristic equation are given by~(\ref{roots:super-c}).  The
solution $u(t)$ of (\ref{pbm:h-eqn})--(\ref{pbm:h-data}) is the sum of
four functions $v_{1}(t)$, $v_{2}(t)$, $w_{1}(t)$, $w_{2}(t)$ whose
components are 
$$\vt_{1}(t,\xi)=
-\ut_{0}(\xi)\cdot\frac{x_{2}(\xi)}{x_{1}(\xi)-x_{2}(\xi)}
\cdot e^{-x_{1}(\xi)t},
\hspace{1.5em}
\vt_{2}(t,\xi)=\ut_{0}(\xi)\cdot\frac{x_{1}(\xi)}{x_{1}(\xi)-x_{2}(\xi)}
\cdot e^{-x_{2}(\xi)t},$$
$$\wt_{1}(t,\xi)=
-\ut_{1}(\xi)\cdot\frac{1}{x_{1}(\xi)-x_{2}(\xi)}\cdot e^{-x_{1}(\xi)t},
\hspace{1.5em}
\wt_{2}(t,\xi)=\ut_{1}(\xi)\cdot\frac{1}{x_{1}(\xi)-x_{2}(\xi)}
\cdot e^{-x_{2}(\xi)t}.$$

The regularity of these four functions follows quite easily from
Lemma~\ref{lemma:main-h}, applied with straightforward choices of the
functions $z_{0}(\xi)$, $c(\xi)$, $g(t,\xi)$.  We skip the elementary
checks, but we sum up the results in the table below.  The first
columns show the values of $P$, $Q$, $R$, $S$ for which the
assumptions are satisfied.  Then in the second column from the right
we write the optimal value of $\alpha$ for which the $m$-th derivative
of the function in that row lies in $C^{0}([0,+\infty),D(A^{\alpha}))$
(according to statement~(2) of Lemma~\ref{lemma:main-h} this value is
$P+Q-mR$, provided of course that it is nonnegative).  Finally, in the
last column we state the values of $\sigma$ for which the
corresponding function is in $C^{\infty}((0,+\infty),D(A^{\alpha}))$
for every $\alpha\geq 0$ (according to statement~(3) of
Lemma~\ref{lemma:main-h} this condition is always equivalent to
$S>0$).

$$\renewcommand{\arraystretch}{1.3}
\begin{array}{|c||c|c|c|c||c|c|}
	\hline
	  & P & Q & R & S & \mbox{$m$-th deriv.\ in $D(A^{\alpha})$} &
	  \mbox{$C^{\infty}$ in $D(A^{\infty})$} \\
	\hline\hline
	\vt_{1}(t,\xi) & \alpha_{0} & 2\sigma-1 & \sigma & \sigma & 
	\alpha_{0}+2\sigma-1 -m\sigma  & \mbox{always}   \\
	\hline
	\vt_{2}(t,\xi) & \alpha_{0} & 0 & 1-\sigma & 1-\sigma & 
	\alpha_{0}-m(1-\sigma) & \mbox{if }\sigma<1  \\
	\hline
	\wt_{1}(t,\xi) & \alpha_{1} & \sigma & \sigma & \sigma & 
	\alpha_{1}+\sigma-m\sigma & \mbox{always}  \\
	\hline
	\wt_{2}(t,\xi) & \alpha_{1} & \sigma & 1-\sigma & 1-\sigma & 
	\alpha_{1}+\sigma-m(1-\sigma) & \mbox{if }\sigma<1  \\
	\hline
\end{array}$$

The regularity of $u(t)$ is the minimal regularity of the four
summands.  Thus we obtain the following results.
\begin{itemize}

	\item  From the second column from the right with $m=0$ we 
	deduce that
	$$u\in
	C^{0}\left([0,+\infty),D(A^{\min\{\alpha_{0},\alpha_{1}+\sigma\}})\right),$$
	while for $m\geq 1$ we deduce that
	$$u^{(m)}\in
	C^{0}\left([0,+\infty),
	D(A^{\min\{\alpha_{0}+\sigma-1,\alpha_{1}\}-(m-1)\sigma})\right).$$
	
	Since $1-\sigma\leq\alpha_{0}-\alpha_{1}\leq\sigma$, this proves
	(\ref{th:h:phsp}) and (\ref{th:h:0-m}) in the supercritical
	regime.

	\item  For $\sigma<1$ all four functions are in 
	$C^{\infty}((0,+\infty),D(A^{\alpha}))$ for every $\alpha\geq 0$, 
	so that the same is true for their sum $u(t)$. This proves 
	(\ref{th:h:>0}) in the supercritical regime.

	\item For $\sigma\geq 1$ we obtain that the $m$-th derivatives of
	$\vt_{2}(t,\xi)$ and $\wt_{2}(t,\xi)$ exist in the space
	$C^{0}([0,+\infty),D(A^{\alpha_{0}+m(\sigma-1)}))$.  Since the
	other two functions are always in
	$C^{\infty}((0,+\infty),D(A^{\alpha}))$ for every $\alpha\geq 0$,
	this proves (\ref{th:h:>0'}).
\end{itemize}

\subparagraph{\textmd{\emph{Critical dissipation}}}

Let us consider the case where $\sigma=1/2$ and $\delta=1$, so that
the roots of the characteristic polynomial are given by~(\ref{roots:c}).
In this case the solution $u(t)$ of
(\ref{pbm:h-eqn})--(\ref{pbm:h-data}) is the sum of three functions
$v_{1}(t)$, $v_{2}(t)$, $w(t)$ whose components are 
$$\vt_{1}(t,\xi)= \ut_{0}(\xi)\cdot 1\cdot e^{-\laxi^{1/2}t}, 
\hspace{3em}
\vt_{2}(t,\xi)=\ut_{0}(\xi)\cdot 1 \cdot \laxi^{1/2}t\,e^{-\laxi^{1/2}t},$$
$$\wt(t,\xi)=\ut_{1}(\xi)\cdot\frac{1}{\laxi^{1/2}}\cdot\laxi^{1/2}t\,e^{-\laxi^{1/2}t}.$$

The regularity of these three functions follows from
Lemma~\ref{lemma:main-h} as before.  We point out that we choose
$g(t,\xi):=\laxi^{1/2}t\,e^{-\laxi^{1/2}t}$ in the case of $\vt_{2}(t,\xi)$
and $\wt(t,\xi)$.  We sum up the results in the table below.
$$\renewcommand{\arraystretch}{1.3}
\begin{array}{|c||c|c|c|c||c|c|}
	\hline
	& P & Q & R & S & \mbox{$m$-th deriv.\ in $D(A^{\alpha})$} &
	\mbox{$C^{\infty}$ in $D(A^{\infty})$} \\
	\hline\hline
	\vt_{1}(t,\xi) & \alpha_{0} & 0 & 1/2 & 1/2 & 
	\alpha_{0}-m/2  & \mbox{always}   \\
	\hline
	\vt_{2}(t,\xi) & \alpha_{0} & 0 & 1/2 & 1/2 & 
	\alpha_{0}-m/2 & \mbox{always}  \\
	\hline
	\wt(t,\xi) & \alpha_{1} & 1/2 & 1/2 & 1/2 & 
	\alpha_{1}+1/2-m/2 & \mbox{always}  \\
	\hline
\end{array}$$

The regularity of $u(t)$ is the minimal regularity of the three
summands.  Looking at the second column from the right, and taking
into account that in this case $\alpha_{0}=\alpha_{1}+1/2$, we obtain
that
$$u^{(m)}\in C^{0}([0,+\infty),D(A^{\alpha_{0}-m/2}))
\quad\quad
\forall m\in\n.$$

This establishes (\ref{th:h:phsp}) and (\ref{th:h:0-m}) in the
critical regime.

Moreover, since $S>0$, all three functions are in
$C^{\infty}((0,+\infty),D(A^{\alpha}))$ for every $\alpha\geq 0$,
which proves (\ref{th:h:>0}) in the critical regime.

\subparagraph{\textmd{\emph{Subcritical dissipation}}}

Let us consider the case where $\sigma<1/2$, or $\sigma=1/2$ and
$\delta<1$, so that for $\lambda$ large enough the roots of the
characteristic polynomial are given by~(\ref{roots:sub-c}).  The
solution $u(t)$ of (\ref{pbm:h-eqn})--(\ref{pbm:h-data}) is the sum of
three functions $v_{1}(t)$, $v_{2}(t)$, $w(t)$ whose components are
$$\vt_{1}(t,\xi)= \ut_{0}(\xi)\cdot 1\cdot e^{-a(\xi)t}\cos(b(\xi)t),
\hspace{3em}
\vt_{2}(t,\xi)=\ut_{0}(\xi)\cdot\frac{a(\xi)}{b(\xi)}\cdot e^{-a(\xi)t}\sin(b(\xi)t),$$
$$\wt(t,\xi)=
\ut_{1}(\xi)\cdot\frac{1}{b(\xi)}\cdot e^{-a(\xi)t}\sin(b(\xi)t).$$

The regularity of these three functions follows from
Lemma~\ref{lemma:main-h} as in the previous cases.  We sum up the
results in the table below. We just point out that the growth of 
derivatives of $g(t,\xi)$, represented by the parameter $R$, is due 
to the terms $b(\xi)$ coming from the trigonometric part. Since 
$b(\xi)\sim\laxi^{1/2}$, we have that $R=1/2$, and this is the reason 
why the derivative gap is always 1/2 in the subcritical regime. 

$$\renewcommand{\arraystretch}{1.3}
\begin{array}{|c||c|c|c|c||c|c|}
	\hline
	& P & Q & R & S & \mbox{$m$-th deriv.\ in $D(A^{\alpha})$} &
	\mbox{$C^{\infty}$ in $D(A^{\infty})$} \\
	\hline\hline
	\vt_{1}(t,\xi) & \alpha_{0} & 0 & 1/2 & \sigma & 
	\alpha_{0}-m/2  & \mbox{if }\sigma>0   \\
	\hline
	\vt_{2}(t,\xi) & \alpha_{0} & 1/2-\sigma & 1/2 & \sigma & 
	\alpha_{0}+1/2-\sigma-m/2 & \mbox{if }\sigma>0  \\
	\hline
	\wt(t,\xi) & \alpha_{1} & 1/2 & 1/2 & \sigma & 
	\alpha_{1}+1/2-m/2 & \mbox{if }\sigma>0  \\
	\hline
\end{array}$$

The regularity of $u(t)$ is the minimal regularity of the three
summands.  Since $\alpha_{0}=\alpha_{1}+1/2$ also in the subcritical
case, as before we obtain that
$$u^{(m)}\in
C^{0}([0,+\infty),D(A^{\alpha_{0}-m/2})) \quad\quad
\forall m\in\n,$$
which proves (\ref{th:h:phsp}) and (\ref{th:h:0-m}) in the subcritical
regime.

Moreover, all three functions are in
$C^{\infty}((0,+\infty),D(A^{\alpha}))$ for every $\alpha\geq 0$
provided that $S>0$, hence when $\sigma>0$.  This proves
(\ref{th:h:>0}) in the subcritical regime.\qed

\section{Proofs in the non-homogeneous case}\label{sec:non-autonomous}

In this section we keep the notation of Section 3. 

\subsection{A  lemma adapted to the forced case} Solutions of
(\ref{eqn:uxi}) with null initial data can always be written as
integrals of $\ft(t,\xi)$ multiplied by some convolution kernel
containing exponential terms.  Let us address the
regularity of these objects.

\begin{lemma}\label{lemma:main}
	Let $H$, $A$, $(M,\mu)$, $\laxi$ be as in
	section~\ref{sec:notation}.  Let us assume that $\laxi\geq 1$ for
	every $\xi\in M$. Let $T>0$, and let $f\in L^{\infty}((0,T),H)$ be a bounded
	function whose components we denote by $\ft(t,\xi)$.  Let
	$y(\xi)$ be a measurable real function on $M$, let $\eta(\xi)$ be 
	a measurable positive function on $M$, and let 
	$\psi(t,\xi)$ be a function defined in $[0,T]\times M$ which is 
	continuous with respect to $t$ and measurable with respect to 
	$\xi$. 
	
	Then the following statements hold true.	
	
	\begin{enumerate}
		\renewcommand{\labelenumi}{(\arabic{enumi})} 
		
		\item Let us assume that there exists a constant $M_{1}$ such 
		that
		\begin{equation}
			|y(\xi)|\cdot|\psi(t,\xi)|\leq M_{1}
			\quad\quad
			\forall(t,\xi)\in[0,T]\times M.
			\label{hp:lemma-main-1}
		\end{equation}
		
		Then the integral
		\begin{equation}
			\zt(t,\xi):=
			y(\xi)\int_{0}^{t}e^{-\eta(\xi)(t-s)}\psi(t-s,\xi)\ft(s,\xi)\,ds
			\label{defn:zt}
		\end{equation}
		defines a function $\zt$ which corresponds, under the usual
		identification of $L^{2}(M,\mu)$ with $H$, to a function $z\in
		C^{0}([0,T],H)$.
		
		\item If in addition to (\ref{hp:lemma-main-1}) we assume 
		also the existence of real numbers
		$\alpha\geq 0$, $b\in[0,1)$, $c\geq 0$, $M_{\alpha,b,c}$ such that
		\begin{equation}
			\laxi^{\alpha}\cdot|y(\xi)|\cdot
			|\psi(t,\xi)|\leq 
			M_{\alpha,b,c}\min\{\eta(\xi)^{b},\eta(\xi)^{c}\}
			\quad\quad
			\forall(t,\xi)\in[0,T]\times M,
		\label{hp:lemma-main-2}
		\end{equation}
		then it turns out that $z\in C^{0}([0,T],D(A^{\alpha}))$. 
		Moreover, there
		exists a constant $K_{b,c}$, depending only on $b$ and $c$,
		such that
		\begin{equation}
			|A^{\alpha}z(t)|\leq K_{b,c}\cdot
			M_{\alpha,b,c}\cdot
			\|f\|_{L^{\infty}((0,T),H)}\cdot
			\int_{0}^{t}\min\left\{\frac{1}{s^{b}},\frac{1}{s^{c}}\right\}\,ds
			\label{th:lemma-main}
		\end{equation}
		for every $t\in[0,T]$, and more generally
		\begin{equation}
			|A^{\alpha}z(t)|\leq K_{b,c}\cdot
			M_{\alpha,b,c}\cdot
			\|f\|_{L^{p}((0,T),H)}\cdot
			\left(\int_{0}^{t}
			\min\left\{\frac{1}{s^{qb}},\frac{1}{s^{qc}}\right\}\,ds
			\right)^{1/q}
			\label{th:lemma-main-p}
		\end{equation}
		for every $t\in[0,T]$, and every pair of real 
		numbers $(p,q)$ such that
		\begin{equation}
			p>1,
			\quad\quad
			q>1,
			\quad\quad
			\frac{1}{p}+\frac{1}{q}=1,
			\quad\quad
			qb<1.
			\label{hp:pq}
		\end{equation}
	
		\item As a special case, let us assume that there exists 
		$M_{2}$ such that
		$$|\psi(t,\xi)|\leq M_{2}
		\quad\quad
		\forall(t,\xi)\in[0,T]\times M,$$
		and that there exist $Q\geq 0$, $S\in\re$, and real numbers 
		$M_{3}$, $M_{4}$, $M_{5}$ such that
		\begin{equation}
			\laxi^{Q}|y(\xi)|\leq M_{3},
			\quad\quad\mbox{and}\quad\quad
			0<M_{4}\leq\frac{\eta(\xi)}{\laxi^{S}}\leq M_{5}
			\label{hp:lemma-qs}
		\end{equation}
		for every $\xi\in M$.
		
		Then $z$ is continuous and bounded in the spaces shown in the 
		following table (it is intended that $0<\ep\leq Q+S$ when 
		needed).
		$$\renewcommand{\arraystretch}{1.2}
		\begin{array}{|c|c|c|}
			\hline
			 & \mbox{$z$ continuous in} & \mbox{$z$ bounded in}  \\
			\hline
			S>0 & D(A^{Q+S-\ep}) & 
			D(A^{Q+S-\ep})  \\
			\hline
			S=0 & D(A^{Q}) & D(A^{Q})  \\
			\hline
			S<0 & D(A^{Q}) & D(A^{Q+S-\ep})  \\
			\hline
		\end{array}$$
		
		Here ``$z$ continuous in $D(A^{\alpha})$'' means that $z\in
		C^{0}([0,T],D(A^{\alpha}))$ and $|A^{\alpha}z(t)|$ satisfies
		an estimate such as (\ref{th:lemma-main}) for some $b<1$ and
		$c\geq 0$,
		while ``$z$ bounded in $D(A^{\alpha})$'' means that in
		addition $c>1$, so that the right-hand side of
		(\ref{th:lemma-main}) is bounded independently of $t$.
	\end{enumerate}
\end{lemma}

\paragraph{\textmd{\textit{Proof}}}

Due to the usual identification of $H$ with $L^{2}(M,\mu)$, it turns 
out that 
\begin{equation}
	\left\|\ft(t,\xi)\right\|_{L^{2}(M,\mu)}=|f(t)|
	\quad\quad
	\mbox{for almost every $t\in(0,T)$.}
	\label{norm-equiv}
\end{equation}

\subparagraph{\textmd{\textit{Statement (1)}}}

The $L^{2}(M,\mu)$ norm of the integral defining $\zt(t,\xi)$ in
(\ref{defn:zt}) is less than or equal to the integral of the norm of
the integrand.  Therefore, from (\ref{hp:lemma-main-1}) it follows
that
\begin{eqnarray*}
	\left\|\zt(t,\xi)\right\|_{L^{2}(M,\mu)} & \leq & 
	\int_{0}^{t}\left\|y(\xi)e^{-\eta(\xi)(t-s)}
	\psi(t-s,\xi)\ft(s,\xi)\right\|_{L^{2}(M,\mu)}\,ds   \\
	\noalign{\vspace{1ex}}
	 & \leq & M_{1}\int_{0}^{t}
	 \left\|\ft(s,\xi)\right\|_{L^{2}(M,\mu)}\,ds  \\
	 \noalign{\vspace{1ex}}
	 & \leq & M_{1}\cdot t\cdot\|f\|_{L^{\infty}((0,T),H)}
\end{eqnarray*}
for every $t\in[0,T]$, which proves that $\zt\in
L^{\infty}\left((0,T),L^{2}(M,\mu)\right)$.  The continuity with
values in the same space follows from the continuity of the map
$$t\to e^{-\eta(\xi)(t-s)}\psi(t-s,\xi)$$
and Lebesgue's theorem.

\subparagraph{\textmd{\textit{Statement (2)}}}

With a simple variable change, we write $\zt(t,\xi)$ in
the form
\begin{equation}
	\zt(t,\xi)=y(\xi)\int_{0}^{t}e^{-\eta(\xi)s}\psi(s,\xi)\ft(t-s,\xi)\,ds.
	\label{stella}
\end{equation}

Then we set
\begin{equation}
	K_{b,c}:=\max\left\{e^{-x}\cdot\max\{x^{b},x^{c}\}:x\geq 0\right\}
	\label{defn:Kbc}
\end{equation}
and, for every $\xi\in M$ and every $0\leq s\leq t\leq T$, we consider the 
function
$$\vft(s,t,\xi):=
\laxi^{\alpha}y(\xi)e^{-\eta(\xi)s}\psi(s,\xi)\ft(t-s,\xi).$$

We claim that for every $t\in(0,T)$ the estimate
\begin{equation}
	\|\vft(s,t,\xi)\|_{L^{2}(M,\mu)}\leq K_{b,c}\cdot M_{\alpha,b,c}\cdot
	\min\left\{\frac{1}{s^{b}},\frac{1}{s^{c}}\right\}\cdot
	|f(t-s)|
	\label{est:varphi}
\end{equation}
holds true for almost every $s\in(0,t)$. Indeed from (\ref{defn:Kbc}) 
it follows that
$$e^{-x}\leq K_{b,c}\min\left\{\frac{1}{x^{b}},\frac{1}{x^{c}}\right\}
\quad\quad
\forall x> 0,$$
hence
$$e^{-\eta(\xi)s}\leq 
K_{b,c}\min\left\{\frac{1}{\eta(\xi)^{b}s^{b}},
\frac{1}{\eta(\xi)^{c}s^{c}}\right\}\leq
K_{b,c}\cdot\frac{1}{\min\{\eta(\xi)^{b},\eta(\xi)^{c}\}}\cdot
\min\left\{\frac{1}{s^{b}},\frac{1}{s^{c}}\right\}$$
for every $s\geq 0$. Thus from (\ref{hp:lemma-main-2}) it follows 
that
\begin{eqnarray*}
	|\vft(s,t,\xi)| & \leq & K_{b,c}\cdot
	\frac{\laxi^{\alpha}|	y(\xi)|}{\min\{\eta(\xi)^{b},\eta(\xi)^{c}\}}\cdot
	\min\left\{\frac{1}{s^{b}},\frac{1}{s^{c}}\right\}\cdot
	|\psi(s,\xi)|\cdot|\ft(t-s,\xi)| \\
	\noalign{\vspace{1ex}}
	 & \leq &  K_{b,c}\cdot
	M_{\alpha,b,c}\cdot
	\min\left\{\frac{1}{s^{b}},\frac{1}{s^{c}}\right\}\cdot
	|\ft(t-s,\xi)|.
\end{eqnarray*}

The coefficient of $|\ft(t-s,\xi)|$ is independent of $\xi$.
Therefore, if we integrate with respect to $\xi$ and we exploit
(\ref{norm-equiv}), we obtain (\ref{est:varphi}).

Now we are ready to prove (\ref{th:lemma-main}).  Indeed from
(\ref{stella}) we obtain that
$$|A^{\alpha}z(t)|= \|\laxi^{\alpha}\zt(t,\xi)\|_{L^{2}(M,\mu)}=
\left\|\int_{0}^{t}\vft(s,t,\xi)\,ds\right\|_{L^{2}(M,\mu)}
\leq\int_{0}^{t}\left\|\vft(s,t,\xi)\right\|_{L^{2}(M,\mu)}\,ds$$
for every $t\in[0,T]$.
Plugging (\ref{est:varphi}) into the last term we easily 
obtain (\ref{th:lemma-main}). 

Let us assume now that $p$ and $q$ satisfy (\ref{hp:pq}).  Then
H\"{o}lder's inequality gives
$$\hspace{-2em}
\int_{0}^{t}\left\|\vft(s,t,\xi)\right\|_{L^{2}(M,\mu)}\,ds \leq
K_{b,c}\cdot M_{\alpha,b,c}\int_{0}^{t}
\min\left\{\frac{1}{s^{b}},\frac{1}{s^{c}}\right\} |f(t-s)|\,ds $$
$$\hspace{2em}
\leq K_{b,c}\cdot M_{\alpha,b,c} \left(\int_{0}^{t}
\min\left\{\frac{1}{s^{qb}},\frac{1}{s^{qc}}\right\}
ds\right)^{1/q}\cdot
\left(\int_{0}^{t}|f(t-s)|^{p}\,ds\right)^{1/p},$$
which easily implies (\ref{th:lemma-main-p}).

For the time being we just proved that
$z\in L^{\infty}((0,T),D(A^{\alpha}))$.  Now we want to prove that $z$
is actually continuous with values in the same space.  To this end, we
show that $A^{\alpha}z(t)$ is the uniform limit of continuous
functions.  Let us set
$$M_{n}:=\{\xi\in M:\laxi\leq n\},$$
and let us set $\chi_{n}(\xi)=1$ if $\xi\in M_{n}$ and 
$\chi_{n}(\xi)=0$ otherwise. It is well-known that $\mu(M\setminus 
M_{n})\to 0$ as $n\to +\infty$. Let us set
\begin{equation}
	\zt_{n}(t,\xi):=\zt(t,\xi)\chi_{n}(\xi),
	\hspace{3em}
	\ft_{n}(t,\xi):=\ft(t,\xi)\chi_{n}(\xi).
	\label{defn:zn-fn}
\end{equation}

Due to the boundedness of $\laxi$ in $M_{n}$, the same argument of
statement~(1) proves that $\laxi^{\alpha}\cdot\zt_{n}(t,\xi)\in
C^{0}\left([0,T],L^{2}(M,\mu)\right)$, which is equivalent to saying
that $\zt_{n}$ corresponds to a function $z_{n}\in
C^{0}([0,T],D(A^{\alpha}))$.  Moreover, from Lebesgue's theorem it
follows that
$$\ft_{n}\to\ft
\quad
\mbox{in }L^{p}((0,T),L^{2}(M,\mu))$$
for every $p\geq 1$ (but not necessarily for $p=+\infty$).  Once
again, this is equivalent to saying that the sequence $\ft_{n}(t,\xi)$
represents the components of a sequence of functions $f_{n}\in
L^{\infty}((0,T),H)$ such that
\begin{equation}
	f_{n}\to f
	\quad
	\mbox{in }L^{p}((0,T),H)
	\label{G->g}
\end{equation}

Let us choose $p$ large enough so that all assumptions
(\ref{hp:pq}) are satisfied.  By the same argument as before, we
obtain
$$\left|A^{\alpha}(z(t)-z_{n}(t))\right|\leq K_{b,c}\cdot
M_{\alpha,b,c} \left(\int_{0}^{t}
\min\left\{\frac{1}{s^{qb}},\frac{1}{s^{qc}}\right\}
ds\right)^{1/q}\cdot 
\|f-f_{n}\|_{L^{p}((0,T),H)}$$
for every $t\in[0,T]$, hence by (\ref{G->g}) we can conclude that
$$\lim_{n\to +\infty}\sup_{t\in[0,T]}
\left|A^{\alpha}(z(t)-z_{n}(t))\right|=0.$$

This is equivalent to saying that $z_{n}(t)\to z(t)$ uniformly in
$C^{0}\left([0,T],D(A^{\alpha})\right)$, and the result follows since
a uniform limit of continuous functions is continuous.

\subparagraph{\textmd{\textit{Statement (3)}}}

All conclusions of this statement follow from
(\ref{th:lemma-main}) with suitable choices of the parameters
$b$ and $c$, which we list below.

\begin{center}
	\renewcommand{\arraystretch}{1.1}
	\begin{tabular}{|c||c|c||c|}
		\hline
		Assumption & $b$ & $c$ & Conclusion  \\
		\hline\hline
		$S>0$ & $1-\ep/S$ & 2 & $z$ continuous and bounded in 
		$D(A^{Q+S-\ep})$  \\
		\hline
		$S=0$ & 0 & 2 & $z$ continuous and bounded in 
		$D(A^{Q})$  \\
		\hline
		$S<0$ & 0 & 0 & $z$ continuous in $D(A^{Q})$  \\
		\hline
		$S<0$ & 0 & $1-\ep/S$ & $z$ bounded in $D(A^{Q+S-\ep})$  \\
		\hline
	\end{tabular}
\end{center}

The verification of (\ref{hp:lemma-main-1}) and
(\ref{hp:lemma-main-2}) in all these cases, with the value of the
exponent $\alpha$ given in the conclusion, is a straightforward check.
We just point out that in the case $S>0$ and $S=0$ one can replace
$c=2$ with any $c>1$, and that in the last line one has that $c>1$
because $S<0$.\qed

\begin{rmk}\label{rmk:lemma-growth}
	\begin{em}
		Lemma~\ref{lemma:main} is stated as a local result in a
		bounded interval $[0,T]$, but it is designed to provide also
		global-in-time estimates when $f$ is defined for all positive
		times.  Let us examine for example the third statement of
		Lemma~\ref{lemma:main}. When $S<0$ and $\alpha\leq Q$, one can
		always apply estimate (\ref{th:lemma-main}) with $b=0$ and
		$c=(\alpha-Q)/S$.  If in addition $\alpha<Q+S$ we obtain
		$c>1$, which gives the boundedness as stated.  If
		$Q+S\leq\alpha\leq Q$ (remember that $S<0$) we obtain that
		$0\leq c\leq 1$.  In this case the right-hand side
		of~(\ref{th:lemma-main}) is not bounded independently of $t$,
		but its growth can be explicitly estimated as follows
		$$
		|A^{\alpha}z(t)|\leq C\cdot
		\|f\|_{L^{\infty}((0,T),H)}\cdot
		\left\{
		\begin{array}{ll}
			t^{(Q+S-\alpha)/S} & 
			\mbox{if }Q+S<\alpha\leq Q , \\
			\noalign{\vspace{0.5ex}}
			\log(1+t) & \mbox{if }\alpha=Q+S,
		\end{array}
		\right.
		$$
		where $C$ is a suitable constant independent of $f$ and $t$.
		This estimate, applied with suitable values of $Q$ and $S$
		(those exploited in the proof of Theorems~\ref{thm:nh-local}
		and~\ref{thm:nh-bound}), leads to~(\ref{th:growth}).
		
	\end{em}
\end{rmk}

\subsection{Proof of Theorem~\ref{thm:nh-local} and 
Theorem~\ref{thm:nh-bound}}

It is easy to express the components $\ut'(t,\xi)$ of the solution in
terms of integrals involving the components $\ft(t,\xi)$ of the forcing
term $f(t)$ and the fundamental solutions of the associated homogenous
equation.  As in the proof of Theorem~\ref{thm:homog}, we
distinguish three cases depending on the discriminant of the
characteristic polynomial~(\ref{char-pol-xi}).

\subparagraph{\textmd{\emph{Supercritical dissipation}}}

Let us consider the case where $\sigma>1/2$, or $\sigma=1/2$ and
$\delta>1$, so that for $\lambda$ large enough the roots of the
characteristic equation are given by~(\ref{roots:super-c}).  The
solution $u(t)$ of (\ref{pbm:nh-eqn})--(\ref{pbm:nh-data}) is the sum
of two functions $v_{1}(t)$ and $v_{2}(t)$ whose components are
$$\vt_{i}(t,\xi)=
\frac{(-1)^{i}}{x_{1}(\xi)-x_{2}(\xi)}
\int_{0}^{t}e^{-x_{i}(\xi)(t-s)}\ft(s,\xi)\,ds
\hspace{3em}
i\in\{1,2\}.$$

The regularity of these two functions follows quite easily from
Lemma~\ref{lemma:main}, applied with straightforward choices of
the functions $y(\xi)$, $\eta(\xi)$, and $\psi(t,\xi)\equiv 1$.  We skip
the elementary checks, but we sum up the results in the table below.
We distinguish the three cases $\sigma<1$, $\sigma=1$, $\sigma>1$, and
for each function we show the values of $Q$ and $S$ for which
assumption (\ref{hp:lemma-qs}) is satisfied, and the conclusions 
according to statement~(3) of Lemma~\ref{lemma:main}.

$$\renewcommand{\arraystretch}{1.3}
\begin{array}{|c||c||c|c|c|c|}
	\hline
	  \mbox{Assumption} & \mbox{Function} & Q & S & \mbox{Continuous
	  in } & \mbox{Bounded in} \\
	  \hline\hline
	  \raisebox{-1.8ex}[-6ex][0ex]{$\sigma<1$} & 
	  v_{1}(t) & \sigma &\sigma & D(A^{2\sigma-\ep}) & D(A^{2\sigma-\ep})   \\
	  \cline{2-6}
	   & 
	  v_{2}(t) & \sigma & 1-\sigma & D(A^{1-\ep}) & D(A^{1-\ep})   \\
	  \hline\hline
	  \raisebox{-1.8ex}[-6ex][0ex]{$\sigma=1$} & 
	  v_{1}(t) & 1 &1 & D(A^{2-\ep}) & D(A^{2-\ep})   \\
	  \cline{2-6}
	   & 
	  v_{2}(t) & 1 & 0 & D(A) & D(A)   \\
	  \hline\hline
	  \raisebox{-1.8ex}[-6ex][0ex]{$\sigma>1$} & 
	  v_{1}(t) & \sigma &\sigma & D(A^{2\sigma-\ep}) & D(A^{2\sigma-\ep})   \\
	  \cline{2-6}
	   & 
	  v_{2}(t) & \sigma & 1-\sigma & D(A^{\sigma}) & D(A^{1-\ep})   \\
	\hline
\end{array}$$

The best space where $u(t)$ is continuous or bounded is always the
maximal space where both $v_{1}(t)$ and $v_{2}(t)$ fulfil the same
property.  This proves the conclusions for $u(t)$ required by
Theorem~\ref{thm:nh-local} and Theorem~\ref{thm:nh-bound} in the
supercritical regime.

Analogously, the derivative $u'(t)$ of the solution is the sum of two
functions $w_{1}(t)$ and $w_{2}(t)$, whose components are
$$\wt_{i}(t,\xi)=
\frac{(-1)^{i+1}x_{i}(\xi)}{x_{1}(\xi)-x_{2}(\xi)}
\int_{0}^{t}e^{-x_{i}(\xi)(t-s)}\ft(s,\xi)\,ds
\hspace{3em}
i\in\{1,2\}.$$

The regularity of $w_{1}(t)$ and $w_{2}(t)$ follows in the same way 
from Lemma~\ref{lemma:main}, as shown in the following table. 

$$\renewcommand{\arraystretch}{1.3}
\begin{array}{|c||c||c|c|c|c|}
	\hline
	  \mbox{Assumption} & \mbox{Function} & Q & S & \mbox{Continuous
	  in } & \mbox{Bounded in} \\
	  \hline\hline
	  \raisebox{-1.8ex}[-6ex][0ex]{$\sigma<1$} & 
	  w_{1}(t) & 0 &\sigma & D(A^{\sigma-\ep}) & D(A^{\sigma-\ep})   \\
	  \cline{2-6}
	   & 
	  w_{2}(t) & 2\sigma-1 & 1-\sigma & D(A^{\sigma-\ep}) & D(A^{\sigma-\ep})   \\
	  \hline\hline
	  \raisebox{-1.8ex}[-6ex][0ex]{$\sigma=1$} & 
	  w_{1}(t) & 0 & 1 & D(A^{1-\ep}) & D(A^{1-\ep})   \\
	  \cline{2-6}
	   & 
	  w_{2}(t) & 1 & 0 & D(A) & D(A)   \\
	  \hline\hline
	  \raisebox{-1.8ex}[-6ex][0ex]{$\sigma>1$} & 
	  w_{1}(t) & 0 &\sigma & D(A^{\sigma-\ep}) & D(A^{\sigma-\ep})   \\
	  \cline{2-6}
	   & 
	  w_{2}(t) & 2\sigma-1 & 1-\sigma & D(A^{2\sigma-1}) & D(A^{\sigma-\ep})   \\
	\hline
\end{array}$$

As before, this is enough to prove the conclusions for $u'(t)$
required by Theorem~\ref{thm:nh-local} and Theorem~\ref{thm:nh-bound}
in the supercritical regime.

\subparagraph{\textmd{\emph{Critical dissipation}}}

Let us consider the case where $\sigma=1/2$ and $\delta=1$, so that
the roots of the characteristic equation are given by~(\ref{roots:c}).
In this case the solution $u(t)$ of
(\ref{pbm:nh-eqn})--(\ref{pbm:nh-data}) has components
$$\ut(t,\xi)=\int_{0}^{t}(t-s)e^{-\laxi^{1/2}(t-s)}\ft(s,\xi)\,ds.$$

The right-hand side can be written in the form of the right-hand side 
of (\ref{defn:zt}) with
$$y(\xi):=\frac{1}{\laxi^{1/2}},
\quad\quad
\eta(\xi):=\frac{\laxi^{1/2}}{2},
\quad\quad
\psi(t,\xi):=\laxi^{1/2}t\exp\left(-\frac{\laxi^{1/2}}{2}t\right).$$

Since assumption (\ref{hp:lemma-qs}) is satisfied with $Q=S=1/2$, we
deduce that $u(t)$ is continuous and bounded with values in
$D(A^{1-\ep})$, as required by Theorem~\ref{thm:nh-local} and
Theorem~\ref{thm:nh-bound} in the critical regime.

As for the derivative $u'(t)$, its components are
$$\ut'(t,\xi)=
-\laxi^{1/2}\int_{0}^{t}(t-s)e^{-\laxi^{1/2}(t-s)}\ft(s,\xi)\,ds
+\int_{0}^{t}e^{-\laxi^{1/2}(t-s)}\ft(s,\xi)\,ds.$$

For the first term, we set $y(\xi)\equiv 1$, and we define
$\eta(\xi)$ and $\psi(t,\xi)$ as before.  For the second term, we set
$y(\xi)=\psi(t,\xi)\equiv 1$ and $\eta(\xi)=\laxi^{1/2}$.  In both
cases assumption (\ref{hp:lemma-qs}) is satisfied with $Q=0$ and
$S=1/2$, and therefore both terms are continuous and bounded in
$D(A^{1/2-\ep})$, as required by Theorem~\ref{thm:nh-local} and
Theorem~\ref{thm:nh-bound} in the critical regime.

\subparagraph{\textmd{\emph{Subcritical dissipation}}}

Let us consider the case where $\sigma\in[0,1/2)$, or $\sigma=1/2$ and
$\delta\in(0,1)$, so that for $\lambda$ large enough the roots of the
characteristic equation are given by~(\ref{roots:sub-c}). The
components of the solution $u(t)$ are
$$\ut(t,\xi)=\frac{1}{b(\xi)}\int_{0}^{t}
e^{-a(\xi)(t-s)}\sin(b(\xi)(t-s))\ft(s,\xi)\,ds.$$

The right-hand side can be written in the form of the right-hand side 
of (\ref{defn:zt}) with
$$y(\xi):=\frac{1}{b(\xi)},
\quad\quad\quad
\eta(\xi):=a(\xi),
\quad\quad\quad
\psi(t,\xi):=\sin(b(\xi)t).$$

It is easy to check that assumption (\ref{hp:lemma-qs}) is satisfied
with $Q=1/2$ and $S=\sigma$.  Thus from Lemma~\ref{lemma:main} it
turns out that $u(t)$ is continuous and bounded in
$D(A^{\sigma+1/2-\ep})$ if $\sigma>0$, and in $D(A^{1/2})$ if
$\sigma=0$, as required by Theorem~\ref{thm:nh-local} and
Theorem~\ref{thm:nh-bound} in the subcritical regime.

As for the derivative $u'(t)$, its components are
\begin{eqnarray*}
	\ut'(t,\xi) & = &
	\int_{0}^{t}e^{-a(\xi)(t-s)}\cos(b(\xi)(t-s))\ft(s,\xi)\,ds \\
	 &  & -\frac{a(\xi)}{b(\xi)}\int_{0}^{t}e^{-a(\xi)(t-s)}
	\sin(b(\xi)(t-s))\ft(s,\xi)\,ds.
\end{eqnarray*}

Once again, we apply Lemma~\ref{lemma:main} to both terms, with
straightforward choices of the parameters.  For $\sigma>0$ we
obtain that the first term is continuous and bounded in
$D(A^{\sigma-\ep})$ (since $Q=0$ and $S=\sigma$), while the second
term is continuous and bounded in $D(A^{1/2-\ep})$ (since
$Q=1/2-\sigma$ and $S=\sigma$).  For $\sigma=0$ we obtain that the
first term is continuous and bounded in $H$ (since $Q=0$ and $S=0$),
while the second term is continuous and bounded in $D(A^{1/2})$ (since
$Q=1/2$ and $S=0$).  Therefore, the regularity and boundedness of
$u'(t)$ is the same of the first term, and it is the same required by
Theorem~\ref{thm:nh-local} and Theorem~\ref{thm:nh-bound} in the
subcritical regime.\qed

\subsection{Proof of Theorem~\ref{thm:nh-limit}}

Let us consider the function
$$E(t):=|A^{\sigma/2}u'(t)|^{2}+|A^{(\sigma+1)/2}u(t)|^{2}.$$

An easy computation shows that
$$E'(t)=-4\delta|A^{\sigma}u'(t)|^{2}+
2\langle A^{\sigma}u'(t),f(t)\rangle\leq
-3\delta|A^{\sigma}u'(t)|^{2}+\frac{1}{\delta}|f(t)|^{2},$$
hence
\begin{equation}
	|A^{\sigma/2}u'(t)|^{2}+|A^{(\sigma+1)/2}u(t)|^{2}
	+3\delta\int_{0}^{t}|A^{\sigma}u'(s)|^{2}\,ds\leq
	\frac{1}{\delta}\int_{0}^{t}|f(s)|^{2}\,ds
	\label{est:E-basic}
\end{equation}
for every $t\in[0,T]$, which proves (\ref{th:nhl:u'2}).

In order to prove (\ref{th:nhl:u2}), we can assume as always that the
operator $A$ is coercive, namely there exists a constant $\nu>0$ such
that $\langle Au,u\rangle\geq\nu|u|^{2}$ for every $u\in D(A)$.  This
allows to estimate $|A^{\alpha}u|$ with $|A^{\beta}u|$ (up to a
constant) whenever $\alpha\leq\beta$.  Now we distinguish two cases.

\paragraph{\textmd{\emph{Case $\sigma\in[1/2,1]$}}}

An easy computation shows that 
\begin{eqnarray*}
	|Au|^{2}+
	\frac{d}{dt}\left(\delta|A^{(\sigma+1)/2}u|^{2}\right) & = &
	|A^{1/2}u'|^{2}-
	\frac{d}{dt}\langle Au,u'\rangle+ 
	\langle f,Au\rangle \\
	 & \leq & |A^{1/2}u'|^{2}-
	 \frac{d}{dt}\langle Au,u'\rangle+
	\frac{1}{2}|Au|^{2}+\frac{1}{2}|f|^{2},
\end{eqnarray*}
hence
$$\frac{1}{2}|Au(t)|^{2}+
\frac{d}{dt}\delta|A^{(\sigma+1)/2}u(t)|^{2}\leq
|A^{1/2}u'(t)|^{2}+\frac{1}{2}|f(t)|^{2}-
\frac{d}{dt}\langle Au(t),u'(t)\rangle.$$

Integrating in $[0,t]$ we obtain that
$$\frac{1}{2}\int_{0}^{t}|Au(s)|^{2}\,ds\leq
\int_{0}^{t}|A^{1/2}u'(s)|^{2}\,ds
+\frac{1}{2}\int_{0}^{t}|f(s)|^{2}\,ds
+|A^{(1-\sigma)/2}u(t)|\cdot|A^{\sigma/2}u'(t)|.$$

Now in the right-hand side we estimate $|A^{1/2}u'(t)|$ with
$|A^{\sigma}u'(t)|$, and $|A^{(1-\sigma)/2}u(t)|$ with
$|A^{(1+\sigma)/2}u(t)|$.  This can be done because $\sigma\geq 1/2$
and the operator $A$ can be assumed to be coercive.  At this point
(\ref{th:nhl:u2}) follows easily from~(\ref{est:E-basic}).

\paragraph{\textmd{\emph{Case $\sigma\in[0,1/2]$}}}

An easy computation shows that
\begin{eqnarray*}
	|A^{\sigma+1/2}u|^{2}+
	\frac{d}{dt}\left(\delta|A^{3\sigma/2}u|^{2}\right) & = &
	|A^{\sigma}u'|^{2}-
	\frac{d}{dt}\langle A^{2\sigma}u,u'\rangle+ 
	\langle f,A^{2\sigma}u\rangle \\
	 & \leq & |A^{\sigma}u'|^{2}-
	 \frac{d}{dt}\langle A^{2\sigma}u,u'\rangle+
	\eta|A^{2\sigma}u|^{2}+\frac{1}{\eta}|f|^{2}.
\end{eqnarray*}

Now we estimate $|A^{2\sigma}u(t)|$ with $|A^{\sigma+1/2}u|$.  This
can be done because $2\sigma\leq\sigma+1/2 $ and the operator $A$ can
be assumed to be coercive.  Thus if we take $\eta$ small enough we
find that 
$$\frac{1}{2}|A^{\sigma+1/2}u(t)|^{2}+
\frac{d}{dt}\left(\delta|A^{3\sigma/2}u(t)|^{2}\right)\leq
|A^{\sigma}u'(t)|^{2}+\frac{1}{\eta}|f(t)|^{2}-
\frac{d}{dt}\langle A^{2\sigma}u(t),u'(t)\rangle.$$

Integrating in $[0,t]$ we obtain that
$$\frac{1}{2}\int_{0}^{t}|A^{\sigma+1/2}u(s)|^{2}\,ds+
\delta|A^{3\sigma/2}u(t)|^{2}$$
$$\leq\int_{0}^{t}|A^{\sigma}u'(s)|^{2}\,ds
+\frac{1}{\eta}\int_{0}^{t}|f(s)|^{2}\,ds
+|\langle A^{2\sigma}u(t),u'(t)\rangle|.$$

Since
$$|\langle A^{2\sigma}u(t),u'(t)\rangle|\leq
|A^{3\sigma/2}u(t)|\cdot|A^{\sigma/2}u'(t)|\leq
\delta|A^{3\sigma/2}u(t)|^{2}+
\frac{1}{4\delta}|A^{\sigma/2}u'(t)|^{2},$$
we easily deduce that
$$\frac{1}{2}\int_{0}^{t}|A^{\sigma+1/2}u(s)|^{2}\,ds\leq
\int_{0}^{t}|A^{\sigma}u'(s)|^{2}\,ds
+\frac{1}{\eta}\int_{0}^{t}|f(s)|^{2}\,ds
+\frac{1}{4\delta}|A^{\sigma/2}u'(t)|^{2}.$$

At this point (\ref{th:nhl:u2}) follows from~(\ref{est:E-basic}).\qed

\subsection{Estimates on the whole line}

The basic tool in the proof of Theorem~\ref{thm:glob-bound} is the
following variant of Lemma~\ref{lemma:main}.

\begin{lemma}\label{lemma:main-infty}
	Let $H$, $A$, $(M,\mu)$, $\laxi$ be as in
	section~\ref{sec:notation}.  Let us assume that $\laxi\geq 1$ for
	every $\xi\in M$.  Let $f\in L^{\infty}(\re,H)$ be a bounded
	function whose components we denote by $\ft(t,\xi)$.  Let $y(\xi)$
	be a measurable real function on $M$, let $\eta(\xi)$ be a
	measurable positive function on $M$, and let $\psi(t,\xi)$ be a
	function defined in $\re\times M$ which is continuous with respect
	to $t$ and measurable with respect to $\xi$.
	
	Then the following statements hold true.	
	
	\begin{enumerate}
		\renewcommand{\labelenumi}{(\arabic{enumi})} 
		
		\item Let us assume that there exist real numbers $\alpha\geq
		0$, $b\in[0,1)$, $c>1$ (this is stronger than the
		corresponding assumption in Lemma~\ref{lemma:main}),
		$M_{\alpha,b,c}$ such that
		$$\laxi^{\alpha}\cdot|y(\xi)|\cdot
		|\psi(t,\xi)|\leq M_{\alpha,b,c}
		\min\left\{\eta(\xi)^{b},\eta(\xi)^{c}\right\}
		\quad\quad
		\forall(t,\xi)\in\re\times M.$$
		
		Then the integral
		$$\zt(t,\xi):=
		y(\xi)\int_{-\infty}^{t}
		e^{-\eta(\xi)(t-s)}\psi(t-s,\xi)\ft(s,\xi)\,ds$$
		defines a function $\zt$ corresponding, under the usual
		identification of $L^{2}(M,\mu)$ with $H$, to a function $z\in
		C^{0}(\re,D(A^{\alpha}))\cap L^{\infty}(\re,D(A^{\alpha}))$.
		
		Moreover, there exists a constant $K_{b,c}$, depending only on
		$b$ and $c$, such that
		\begin{equation}
			|A^{\alpha}z(t)|\leq K_{b,c}\cdot
			M_{\alpha,b,c}\cdot
			\|f\|_{L^{\infty}(\re,H)}\cdot
			\int_{0}^{+\infty}
			\min\left\{\frac{1}{s^{b}},\frac{1}{s^{c}}\right\}\,ds
			\quad\quad
			\forall t\in\re.
			\label{th:lemma-main-infty}
		\end{equation}
		
		\item As a special case, let us assume that
		$\sup\left\{|\psi(t,\xi)|:(t,\xi)\in\re\times M\right\}
		<+\infty$, and that there exist $Q\geq 0$, $S\in\re$, and real
		numbers $M_{3}$, $M_{4}$, $M_{5}$ such that
		(\ref{hp:lemma-qs}) holds true for every $\xi\in M$.
		
		Then $z$ is continuous and bounded in the spaces shown in the 
		following table (it is intended that $0<\ep\leq Q+S$ when 
		needed).
		$$\renewcommand{\arraystretch}{1.2}
		\begin{array}{|c|c|}
			\hline
			 & \mbox{$z$ continuous and bounded in}   \\
			\hline
			S>0 & D(A^{Q+S-\ep})  \\
			\hline
			S=0 & D(A^{Q})  \\
			\hline
			S<0 & D(A^{Q+S-\ep})  \\
			\hline
		\end{array}$$
		
		\item Let us assume in addition that $\psi(t,\xi)\equiv 1$,
		that $f(t)$ is periodic, and in (\ref{hp:lemma-qs}) we have
		that $S<0$.
		
		Then $z(t)$ is continuous and bounded also in the limit space
		$D(A^{Q+S})$.
	\end{enumerate}
		
\end{lemma}

\paragraph{\textmd{\textit{Proof}}}

First of all, with a variable change we rewrite $\zt(t,\xi)$ as
$$\zt(t,\xi)=y(\xi)\int_{0}^{+\infty}e^{-\eta(\xi)s}
\psi(s,\xi)\ft(t-s,\xi)\,ds.$$

Now we are ready to prove our conclusions.  The proof of
(\ref{th:lemma-main-infty}) is analogous to the proof of
(\ref{th:lemma-main}), the only difference being that now the integral
is over $(0,+\infty)$ instead of $(0,t)$.  This gives the boundedness
of $z(t)$.  In order to prove the continuity, we define $z_{n}(t)$ and
$f_{n}(t)$ as in (\ref{defn:zn-fn}), and we prove that $z_{n}(t)$
uniformly converges to $z(t)$ on every closed interval $[A,B]$.

Arguing as in the proof of Lemma~\ref{lemma:main} we obtain that
\begin{equation}
	\left|A^{\alpha}(z(t)-z_{n}(t))\right|\leq K_{b,c}\cdot
	M_{\alpha,b,c}\int_{0}^{+\infty}
	\min\left\{\frac{1}{s^{b}},\frac{1}{s^{c}}\right\}
	|f(t-s)-f_{n}(t-s)|\,ds.
	\label{est:A-alpha-infty}
\end{equation}

Let us fix now any $\ep>0$, and let us split the integral in the
right-hand side of (\ref{est:A-alpha-infty}) as an integral in some
bounded interval $[0,T]$ and the integral in $[T,+\infty)$.  For every
$T\geq 0$ we have that
$$\int_{T}^{+\infty}
\min\left\{\frac{1}{s^{b}},\frac{1}{s^{c}}\right\}
|f(t-s)-f_{n}(t-s)|\,ds\leq
\|f\|_{L^{\infty}(\re,H)}\int_{T}^{+\infty}
\min\left\{\frac{1}{s^{b}},\frac{1}{s^{c}}\right\}\,ds.$$

Therefore, if $T$ is large enough we have that
\begin{equation}
	K_{b,c}\cdot M_{\alpha,b,c}\int_{T}^{+\infty}
	\min\left\{\frac{1}{s^{b}},\frac{1}{s^{c}}\right\}
	|f(t-s)-f_{n}(t-s)|\,ds\leq\frac{\ep}{2}
	\quad\quad
	\forall n\in\n.
	\label{est:e/2-1}
\end{equation}

Let us consider now the integral in $[0,T]$. Let us choose $p$ and $q$ in 
such a way that all assumptions in (\ref{hp:pq}) are satisfied. When 
$t\in[A,B]$ and $s\in[0,T]$ we have that $t-s\in[A-T,B]$, hence
\begin{eqnarray}
	\lefteqn{\hspace{-2em}\int_{0}^{T}
	\min\left\{\frac{1}{s^{b}},\frac{1}{s^{c}}\right\}
	|f(t-s)-f_{n}(t-s)|\,ds}
	\nonumber  \\
	 & \leq & \left(\int_{0}^{T}
	 \min\left\{\frac{1}{s^{qb}},\frac{1}{s^{qc}}\right\}\,ds\right)^{1/q}
	 \left(\int_{0}^{T} |f(t-s)-f_{n}(t-s)|^{p}\,ds\right)^{1/p}
	\nonumber  \\
	 & \leq & \left(\int_{0}^{T}
	 \min\left\{\frac{1}{s^{qb}},\frac{1}{s^{qc}}\right\}\,ds\right)^{1/q}
	 \|f-f_{n}\|_{L^{p}((A-T,B),H)}.
	\nonumber
\end{eqnarray}

Since $f_{n}\to f$ in $L^{p}((A-T,B),H)$, when $n$ is large enough we 
have that
\begin{equation}
	K_{b,c}\cdot M_{\alpha,b,c}
	\int_{0}^{T} \min\left\{\frac{1}{s^{b}},\frac{1}{s^{c}}\right\}
	|f(t-s)-f_{n}(t-s)|\,ds\leq\frac{\ep}{2}.
	\label{est:e/2-2}
\end{equation}

Plugging (\ref{est:e/2-1}) and (\ref{est:e/2-2}) into
(\ref{est:A-alpha-infty}) we obtain that
$$\sup_{t\in[A,B]}\left|A^{\alpha}(z(t)-z_{n}(t))\right|\leq\ep$$
provided that $n$ is large enough. This completes the proof of the 
first statement.

The proof of the second statement is analogous to the proof of 
statement~(3) of Lemma~\ref{lemma:main}.  The unique
difference is that in the case $S<0$ we cannot choose $c=0$, because
now we need $c>1$.  Thus we choose $c=1-\ep/S$ (which is larger than 1
because $S<0$), and we obtain both continuity and boundedness in the
same space, namely $D(A^{Q+S-\ep})$.

It remains to prove statement~(3). If $\psi(t,\xi)\equiv 1$, and 
$f(t)$ is $T_{0}$-periodic for some $T_{0}>0$, then 
$$\zt(t,\xi)=y(\xi)\int_{0}^{+\infty}e^{-\eta(\xi)s}\ft(t-s,\xi)\,ds=
y(\xi)\sum_{n=0}^{\infty}
\int_{nT_{0}}^{(n+1)T_{0}}e^{-\eta(\xi)s}\ft(t-s,\xi)\,ds=$$
$$=
y(\xi)\sum_{n=0}^{\infty}e^{-\eta(\xi)nT_{0}}
\int_{0}^{T_{0}}e^{-\eta(\xi)s}\ft(t-s,\xi)\,ds=
\frac{y(\xi)}{1-e^{-\eta(\xi)T_{0}}}
\int_{0}^{T_{0}}e^{-\eta(\xi)s}\ft(t-s,\xi)\,ds,$$
hence
\begin{eqnarray}
	\laxi^{Q+S}|\zt(t,\xi)| & \leq &
	\frac{\laxi^{Q+S}|y(\xi)|}{1-e^{-\eta(\xi)T_{0}}}
	\int_{0}^{T_{0}}e^{-\eta(\xi)s}|\ft(t-s,\xi)|\,ds 
	\nonumber \\
	\noalign{\vspace{1ex}}
	 & \leq & \frac{\laxi^{Q+S}|y(\xi)|}{1-e^{-\eta(\xi)T_{0}}}
	\int_{0}^{T_{0}}|\ft(t-s,\xi)|\,ds.
	\nonumber
\end{eqnarray}

Since $S<0$, from (\ref{hp:lemma-qs}) we deduce that $\eta(\xi)\leq
M_{6}$ for a suitable constant $M_{6}$, hence there exists a constant
$M_{7}$ such that 
$$\frac{\eta(\xi)}{1-e^{-\eta(\xi)T_{0}}}\leq M_{7}
\quad\quad
\forall\xi\in M.$$

Taking once again (\ref{hp:lemma-qs}) into account, we deduce that 
there exists a constant $M_{8}$ such that
$$\frac{\laxi^{Q+S}|y(\xi)|}{1-e^{-\eta(\xi)T_{0}}}=
\laxi^{Q}|y(\xi)|\cdot
\frac{\laxi^{S}}{\eta(\xi)}\cdot
\frac{\eta(\xi)}{1-e^{-\eta(\xi)T_{0}}}\leq M_{8}
\quad\quad
\forall\xi\in M.$$

It follows that 
$$|A^{Q+S}z(t)|\leq
M_{8}\int_{0}^{T_{0}}|f(t-s)|\,ds\leq
M_{8}\|f\|_{L^{\infty}((0,T_{0}),H)}T_{0} 
\quad\quad
\forall t\in\re,$$
which proves that $z$ is bounded in $D(A^{Q+S})$. The continuity 
follows from the uniform convergence of the sequence $z_{n}(t)$, 
which can be proved exactly in the same way.
\qed

\subsection{Proof of Theorem \ref{thm:glob-bound}}

First we construct a bounded solution on $\re$ fulfilling the various
boundedness and continuity requirements.  The formula defining the
unique bounded solution for a forced exponentially damped system is
our guide since in the $(M, \mu)$ formulation, the restriction of the
system on states with both components supported in $M_{n}:=\{\xi\in
M:\laxi\leq n\}$ is exponentially damped, for every integer $n$.  Let
$\ft(t,\xi)$ denote the components of the forcing term $f(t)$.  Let us
define $\ut(t,\xi)$ as in the proof of Theorem~\ref{thm:nh-local} and
Theorem~\ref{thm:nh-bound}, the only difference being that now the
integration is over $(-\infty,t)$ instead of $(0,t)$.  For example, in
the case of supercritical dissipation, $\ut(t,\xi)$ is given by
$$-\frac{1}{x_{1}(\xi)-x_{2}(\xi)}\int_{-\infty}^{t}
e^{-x_{1}(\xi)(t-s)}\ft(s,\xi)\,ds
+\frac{1}{x_{1}(\xi)-x_{2}(\xi)}\int_{-\infty}^{t}
e^{-x_{2}(\xi)(t-s)}\ft(s,\xi)\,ds.$$

Due to Lemma~\ref{lemma:main-infty}, these two integrals define two
functions (whose sum is a solution to (\ref{pbm:eqn})), whose
regularity is given by statement~(2) of the same lemma.  The
conclusion is that the solution is continuous and bounded in the same
spaces where it was bounded in the case of Theorem~\ref{thm:nh-bound}
(in other words, now there is no difference between the spaces where
the solution is guaranteed to be bounded and the spaces where the
solution is guaranteed to be continuous). The same arguments apply to the critical and subcritical dissipation, 
and to the regularity of derivatives. This completes the proof of the 
first statement. 
\bigskip 

Concerning uniqueness, it is a simple consequence of the fact that the
only solution $v\in C_b^{0}(\re, D(A^{1/2}))\cap C^1_b(\re, H)$ of the
homogeneous equation is $0$.  This follows from the fact that the
restriction of the system on states with both components supported in
$M_{n}:=\{\xi\in M:\laxi\leq n\}$ is exponentially damped, a property
implying that the projections of $v$ on all $M_{n}$ are trivial.  It
follows immediately from uniqueness of the bounded solution that if
$f$ is periodic, the bounded solution is periodic with the same
period.

\bigskip 
It follows obviously from the previous estimates that
uniqueness of the bounded solution $u$ is reinforced by the fact that
the norm of $u$ in $C_b^{0}(\re, D(A^{1/2}))\cap C^1_b(\re, H)$ is bounded
by a constant times the norm of $f$ in $L^\infty(\re, H).$ As a
consequence, if $f$ is almost periodic with values in $H$, the bounded
solution$(u, u')$ is almost periodic from $ \re$ to $D(A^{1/2})\times
H$ with $\exp(u)\subset \exp(f)$.  The property of almost periodicity,
for the same reason, is in fact valid also with values in all the
product spaces where boundedness and continuity has been proved for
the first statement.  \bigskip

It remains to prove that in the periodic case with $\sigma>1$ the
solution is continuous and bounded also in the limit space $D(A)$.  In
this case $u(t)$ is the sum of two terms, whose components are written
above.  The first one is continuous and bounded in
$D(A^{2\sigma-\ep})$ (same proof as in Theorem~\ref{thm:nh-bound}),
and $2\sigma-\ep\geq 1$ if $\ep$ is small enough.  The second term
fits in the framework of statement~(3) of 
Lemma~\ref{lemma:main-infty} with $Q=\sigma$ and $S=1-\sigma<0$. Thus 
we obtain that $u(t)$ is continuous and bounded in $D(A^{Q+S})$, 
which is exactly $D(A)$.\qed

\setcounter{equation}{0}
\section{Counterexamples}\label{sec:counterexamples}

In this section we exhibit all the counterexamples needed in the 
proof of Theorem~\ref{thm:optimal}. To begin with, let us fix some notations.  Let $\{\lk\}$ be an
unbounded sequence of positive eigenvalues of $A$, which we can always
assume to be increasing.  Let $\{e_{k}\}$ be a sequence of
corresponding eigenvectors, which we can always take with unit norm.
Up to restricting to the smallest closed vector subspace containing
the sequence $\{e_{k}\}$, we can assume that $\{e_{k}\}$ is an
orthonormal basis of $H$.  This allows to identify any vector $v\in H$
with the sequence $\{v_{k}\}$ of its components with respect to
$\{e_{k}\}$.  Under this identification, $v\in D(A^{\alpha})$ if and
only if 
$$\sum_{k=0}^{\infty}\lk^{2\alpha}v_{k}^{2}<+\infty.$$

Different statements of Theorem~\ref{thm:optimal} require different
strategies, which now we briefly introduce.

The easiest one is the proof of~(\ref{th:4a}), which is the only 
case where we obtain a solution which lacks in regularity for all 
positive times. In this case a constant forcing term is enough to 
produce the required regularity loss.

Then we pass to examples where we prove lack of regularity on a given
sequence $\{t_{n}\}$ (possibly dense) of positive times.  The
underlying strategy is the same, and it consists in the following two
main steps.
\begin{itemize}
	\item  In the first step we produce a solution with the required 
	regularity loss at a given single time $T>0$. The main point is 
	that this can be done using a forcing term with norm as small as 
	we want, and concentrated on the subspace of $H$ generated by a 
	given countable subset of the eigenvectors $\{e_{n}\}$.

	\item In the second step we begin by partitioning $\{e_{n}\}$ into
	countably many disjoint countable subsets, which thus generate a
	countable set of pairwise orthogonal subsets $H_{n}$ of $H$.  Then
	for each $n$ we apply the result of the first step in order to
	obtain a forcing term $f_{n}(t)$, with values in $H_{n},$ which
	produces a solution with the required regularity loss at time
	$T=t_{n}$.  Since we can take the norm of $f_{n}(t)$ as small as
	we want, we can arrange things so that the series with general
	term $f_{n}(t)$ converges.  The sum $f(t)$ is the required
	forcing term giving rise to a solution $u(t)$ which lacks in
	regularity at all times of the sequence $\{t_{n}\}$.
\end{itemize}

Thanks to this strategy, the proof of statements~(1) and~(2), and of
part~(\ref{th:4b}) of statement~(3), is reduced to the verification of
the first step.  In all cases, this in turn requires a forcing term
defined as the sum of a suitable series of forcing terms, whose
construction and convergence differs from case to case.  When
$\sigma=0$ the series converges because the series of the norms
converges.  In the other two cases, the series converges because its
terms have norms tending to 0 and disjoint supports, namely the
different components are ``activated'' one by one.

Finally, we produce the unbounded solution required by statement~(4).
Once again, the strategy is twofold.  We show that, provided that $t$
is large enough, $|Au(t)|$ can be made larger than a given constant,
even if the forcing term is smaller than a given constant.  Then we
conclude with an argument similar to the second step described above.

All proofs begin with careful asymptotic estimates on solutions to 
the family of ordinary differential equations
\begin{equation}
	\ul''(t)+2\delta\lambda^{\sigma}\ul'(t)+\lambda\ul(t)=\fl(t)
	\quad\quad
	\forall t\geq 0,
	\label{ODE-eqn}
\end{equation}
with a suitable family of forcing terms $\{\fl(t)\}$, and null 
initial data
\begin{equation}
	\ul(0)=\ul'(0)=0.
	\label{ODE-data}
\end{equation}

In turn, the asymptotic behavior of these solutions depends on the
asymptotic behavior of the roots of the characteristic polynomial
(\ref{char-pol}).  In analogy with section~\ref{sec:roots}, these
roots are denoted by $-x_{1,\lambda}$ and $-x_{2,\lambda}$.  Their
expressions and asymptotic behavior are the same stated in
section~\ref{sec:roots}, just with $\lambda$ instead of
$\lambda(\xi)$.  In particular, if $\lambda$ is large enough, the
picture is the following.
\begin{itemize}
	\item In the subcritical regime the characteristic roots are
	complex conjugate numbers of the form $-x_{1,\lambda}=-\al+i\bl$
	and $-x_{2,\lambda}=-\al-i\bl$, with
	\begin{equation}
		\al:=\delta\lambda^{\sigma},
		\hspace{3em}
		\bl:=(\lambda-\delta^{2}\lambda^{2\sigma})^{1/2},
		\label{defn:al-bl-c}
	\end{equation}
	and as a consequence
	\begin{equation}
		\lim_{\lambda\to +\infty}\frac{\bl}{\lambda^{1/2}}=
		\left\{
		\begin{array}{ll}
			1 & \mbox{if $\sigma<1/2$,}  \\
			(1-\delta^{2})^{1/2} & \mbox{if $\sigma=1/2$.}
		\end{array}
		\right.
		\label{th:bl-asympt}
	\end{equation}

	\item In the critical regime it turns out that
	$-x_{1,\lambda}=-x_{2,\lambda}=-\lambda^{1/2}$.

	\item In the supercritical regime the characteristic roots are
	real numbers with
	\begin{equation}
		\lim_{\lambda\to +\infty}\frac{x_{1,\lambda}}{\lambda^{\sigma}}=
		\lim_{\lambda\to +\infty}\frac{\lambda^{1-\sigma}}{x_{2,\lambda}}=
		\left\{
		\begin{array}{ll}
			2\delta & \mbox{if $\sigma>1/2$,}  \\
			\delta+(\delta^{2}-1)^{1/2} & \mbox{if $\sigma=1/2$.}
		\end{array}
		\right.
		\label{th:super-asympt}
	\end{equation}
	
\end{itemize}

\subsection{Proof of~(\ref{th:4a}) in statement~(3)}

Let us consider problem (\ref{ODE-eqn})--(\ref{ODE-data}) with
$\fl(t)\equiv 1$.  Since we are in the case $\sigma\geq 1$, when
$\lambda$ is large enough the solution turns out to be
\begin{equation}
	\ul(t)=\frac{1}{\lambda}+\frac{1}{\lambda}\cdot
	\frac{1}{x_{1,\lambda}-x_{2,\lambda}}\left(
	x_{2,\lambda}e^{-x_{1,\lambda}t}-x_{1,\lambda}e^{-x_{2,\lambda}t}
	\right).
	\label{formula:ul-sc}
\end{equation}

Thanks to (\ref{th:super-asympt}), it is not difficult to show that
\begin{equation}
	\lim_{\lambda\to +\infty}\lambda^{\sigma}\ul(t)=\left\{
	\begin{array}{ll}
		t/(2\delta) & \mbox{if }\sigma>1  \\
		\noalign{\vspace{0.5ex}}
		1-e^{-t/(2\delta)} & \mbox{if }\sigma=1
	\end{array}
	\right.
	\label{ul-lim-4a}
\end{equation}
for all $t\geq 0$.  In both cases, the limit is finite and different
from 0 when $t>0$.

Now let us choose a sequence $a_{k}$ such that
\begin{equation}
	\sum_{k=0}^{\infty}a_{k}^{2}<+\infty
	\label{defn:ak-conv}
\end{equation}
and
\begin{equation}
	\sum_{k=0}^{\infty}\lk^{2\ep}a_{k}^{2}=+\infty
	\quad\quad
	\forall\ep>0.
	\label{defn:ak-div}
\end{equation}

Let us consider the constant forcing term
$$f(t):=\sum_{k=0}^{\infty}a_{k}e_{k}.$$

The series converges because of (\ref{defn:ak-conv}). The 
corresponding solution of  (\ref{pbm:nh-eqn})--(\ref{pbm:nh-data}) is 
$$u(t)=\sum_{k=0}^{\infty}a_{k}u_{\lambda_{k}}(t)e_{k},$$
and in particular
$$\left|A^{\sigma+\ep}u(t)\right|^{2}=
\sum_{k=0}^{\infty}\lk^{2\sigma+2\ep}a_{k}^{2}|u_{\lambda_{k}}(t)|^{2}.$$

Due to (\ref{ul-lim-4a}), this series is equivalent to the series in 
(\ref{defn:ak-div}) for every $t>0$, hence it is divergent for every 
$t>0$ and every $\ep>0$. This proves (\ref{th:4a}).

\subsection{Proof of statement~(1)}\label{ctrex:stat-1}

\paragraph{\textmd{\emph{ODE estimates}}}

In the subcritical case the roots of the characteristic polynomial are
complex conjugate numbers $-\al\pm i\bl$, at least when $\lambda$ is
large enough.  For every such $\lambda$ and every $T>0$, we consider
problem (\ref{ODE-eqn})--(\ref{ODE-data}) with
\begin{equation}
	\fl(t):=\cos\left(\bl(T-t)-\frac{\pi}{4}\right)=
	\frac{\sqrt{2}}{2}\left[\strut
	\cos\left(\bl(T-t)\right)+\sin\left(\bl(T-t)\right)\right].
	\label{defn:fl}
\end{equation}

We claim that
\begin{equation}
	\lim_{\lambda\to +\infty}\lambda^{1/2}\ul(T)=
	\lim_{\lambda\to +\infty}\ul'(T)=
	\frac{\sqrt{2}}{4}\cdot\frac{1}{\delta}\left(1-e^{-\delta T}\right).
	\label{ul-lim-1}
\end{equation}

Indeed the solution of (\ref{ODE-eqn})--(\ref{ODE-data}) is
\begin{equation}
	\ul(t)=\frac{1}{\bl}\int_{0}^{t}
	e^{-\al(t-s)}\sin(\bl(t-s))\fl(s)\,ds,
	\label{soln:ul}
\end{equation}
and its derivative $\ul'(t)$ is 
\begin{equation}
	\ul'(t)=-\al\ul(t)+\int_{0}^{t}
	e^{-\al(t-s)}\cos(\bl(t-s))\fl(s)\,ds.
	\label{soln:ul'}
\end{equation}

In the special case where $\sigma=0$, and $\fl(t)$ is given by 
(\ref{defn:fl}), with the variable change $x=T-s$ we obtain 
that
\begin{equation}
	\lambda^{1/2}\ul(T)=\frac{\sqrt{2}}{2}\frac{\lambda^{1/2}}{\bl}
	\int_{0}^{T}e^{-\delta x}\left(
	\sin^{2}(\bl x)+\sin(\bl x)\cos(\bl x)\right)
	\,dx,
	\label{lim:0-u}
\end{equation}
\begin{equation}
	\ul'(T)=-\delta\ul(T)+\frac{\sqrt{2}}{2}
	\int_{0}^{T}e^{-\delta x}\left(
	\cos^{2}(\bl x)+\sin(\bl x)\cos(\bl x)\right)
	\,dx.
	\label{lim:0-u'}
\end{equation}

Now we have to compute the limits as $\lambda\to +\infty$. The coefficient 
$\lambda^{1/2}\cdot \bl^{-1}$ tends to 1 because of 
(\ref{th:bl-asympt}). As for the integrals, the coefficient $\bl$ in 
the trigonometric terms tends to $+\infty$. This produces a 
homogenization effect, so that
$$\lim_{\lambda\to+\infty}\int_{0}^{T}e^{-\delta x}\sin^{2}(\bl x)\,dx
=\lim_{\lambda\to+\infty}\int_{0}^{T}e^{-\delta x}\cos^{2}(\bl 
x)\,dx=
\frac{1}{2}\int_{0}^{T}e^{-\delta x}\,dx,$$
$$\lim_{\lambda\to+\infty}
\int_{0}^{T}e^{-\delta x}\sin(\bl x)\cos(\bl x)\,dx=0.$$

Plugging these limits into (\ref{lim:0-u}) and (\ref{lim:0-u'}), we
obtain (\ref{ul-lim-1}).

\paragraph{\textmd{\emph{Lack of regularity for a given positive time}}}

Let $T>0$, let $\eta>0$, let $\nu_{k}\to +\infty$ be any unbounded
sequence of eigenvalues of $A$, let $\{\widehat{e}_{k}\}$ be a
corresponding sequence of orthonormal eigenvectors, and let
$\widehat{H}$ be the subspace of $H$ generated by
$\{\widehat{e}_{k}\}$.  We claim that there exists a function $f\in
C^{0}_{b}([0,+\infty),H)$ such that
\begin{equation}
	|f(t)|\leq\eta
	\quad\quad
	\forall t\geq 0,
	\label{est:ft-eta}
\end{equation}
\begin{equation}
	f(t)\in\widehat{H}
	\quad\quad
	\forall t\geq 0,
	\label{est:ft-H}
\end{equation}
and such that the corresponding solution $u$ of
(\ref{pbm:nh-eqn})--(\ref{pbm:nh-data}) satisfies
\begin{equation}
	u(T)\not\in D(A^{1/2+\ep})
	\quad \mbox{and}\quad 
	u'(T)\not\in D(A^{\ep})
	\quad\quad
	\forall \ep>0.
	\label{claim:0-T}
\end{equation}

Indeed let us choose a sequence $a_{k}$ such that
\begin{equation}
	\sum_{k=0}^{\infty}a_{k}^{2}\leq\eta
	\label{defn:ak-conv-eta}
\end{equation}
and
\begin{equation}
	\sum_{k=0}^{\infty}\nu_{k}^{2\ep}a_{k}^{2}=+\infty
	\quad\quad
	\forall\ep>0.
	\label{defn:ak-div-nu}
\end{equation}

Let us consider the family of forcing terms $\fl(t)$ defined in
(\ref{defn:fl}), and the corresponding solutions $\ul(t)$ of problem
(\ref{ODE-eqn})--(\ref{ODE-data}).  Let us set
$$f(t):=\sum_{k=0}^{\infty}a_{k}f_{\nu_{k}}(t)\widehat{e}_{k}.$$

Due to (\ref{defn:ak-conv-eta}), the series converges to a continuous
function $f:[0,+\infty)\to H$ satisfying both (\ref{est:ft-eta}) and
(\ref{est:ft-H}).  The corresponding solution of
(\ref{pbm:nh-eqn})--(\ref{pbm:nh-data}) is of course 
$$u(t)=\sum_{k=0}^{\infty}a_{k}u_{\nu_{k}}(t)\widehat{e}_{k},$$
and in particular
$$\left|A^{1/2+\ep}u(T)\right|^{2}=
\sum_{k=0}^{\infty}\nu_{k}^{1+2\ep}a_{k}^{2}|u_{\nu_{k}}(T)|^{2},$$
$$\left|A^{\ep}u'(T)\right|^{2}=
\sum_{k=0}^{\infty}\nu_{k}^{2\ep}a_{k}^{2}|u_{\nu_{k}}'(T)|^{2}.$$

Due to (\ref{ul-lim-1}), both series are equivalent to the series in
(\ref{defn:ak-div-nu}), hence they diverge for every $\ep>0$.  This
proves (\ref{claim:0-T}).

\paragraph{\textmd{\emph{Lack of regularity for a given sequence of times}}}

We are now ready to prove the conclusions of statement~(1) of
Theorem~\ref{thm:optimal}.  To this end, we partition the given
unbounded sequence $\{\lk\}$ of eigenvalues of $A$ into countably many
disjoint (unbounded) subsequences.  For example, the $n$-th
subsequence could be that of the form $\lambda_{2^{n}(2k+1)}$.

Let $\{t_{n}\}\subseteq(0,+\infty)$ be any sequence of positive times.
For every $n\in\n$, we apply the construction of the previous
paragraph with $T:=t_{n}$, $\eta:=2^{-n}$, and $\{\nu_{k}\}$ equal to
the $n$-th subsequence of $\{\lk\}$.  We call $H_{n}$ the subspace
generated by the corresponding eigenvectors of $A$.  

We obtain an external force $f_{n}\in C^{0}_{b}([0,+\infty),H)$ such
that $f_{n}(t)\in H_{n}$ and $|f_{n}(t)|\leq 2^{-n}$ for every $t\geq
0$, and such that the corresponding solution $u_{n}(t)$ of
(\ref{pbm:nh-eqn})--(\ref{pbm:nh-data}) takes its values in $H_{n}$
and satisfies
\begin{equation}
	u_{n}(t_{n})\not\in D(A^{1/2+\ep})
	\quad \mbox{and}\quad 
	u_{n}'(t_{n})\not\in D(A^{\ep})
	\quad\quad
	\forall \ep>0.
	\label{claim:0-tn}
\end{equation}

Finally, we define 
\begin{equation}
	f(t):=\sum_{k=0}^{\infty}f_{n}(t), 
	\hspace{4em}
	u(t):=\sum_{k=0}^{\infty}u_{n}(t).
	\label{defn:f-u}
\end{equation}

It is easy to see that the first series converges to a bounded continuous 
function $f(t)$, and that $u(t)$ is the corresponding solution of 
(\ref{pbm:nh-eqn})--(\ref{pbm:nh-data}).

Since the spaces $H_{n}$ are pairwise orthogonal, $u(t)$ cannot be
more regular than its projections $u_{n}(t)$ into $H_{n}$, and
therefore both (\ref{th:1a}) and (\ref{th:1b}) follow from
(\ref{claim:0-tn}).

% \clearpage

\subsection{Proof of statement~(2)}

\paragraph{\textmd{\emph{Blow-up triples}}}

We say that a triple $(\sigma,\sigma_{0},\sigma_{1})$ of positive real
numbers satisfies the blow-up condition if there exists families
$\{\tl\}\subseteq(0,+\infty)$ and $\{\fl\}\subseteq
C^{0}_{b}([0,+\infty),H)$ such that
\begin{equation}
	|\fl(t)|\leq 1
	\quad\quad
	\forall t\geq 0,\quad\forall\lambda\geq 0,
	\label{defn:buc-fl}
\end{equation}
\begin{equation}
	\lim_{\lambda\to +\infty}\tl=0,
	\label{defn:buc-tl}
\end{equation}
and the corresponding solutions $\ul(t)$ of
(\ref{ODE-eqn})--(\ref{ODE-data}) satisfy
\begin{equation}
	\lim_{\lambda\to +\infty}
	\lambda^{\sigma_{0}}|\ul(\tl)|\in(0,+\infty)
	\quad\quad\mbox{and}\quad\quad
	\lim_{\lambda\to +\infty}
	\lambda^{\sigma_{1}}|\ul'(\tl)|\in(0,+\infty).
	\label{defn:buc-limits}
\end{equation}

We claim that the triples
$\left(\sigma,\min\left\{\sigma+1/2,1\right\},\sigma\right)$ satisfy
the blow-up condition for every $\sigma\in(0,1)$ .  The verification of
this fact requires several cases.

In the case of supercritical dissipation, namely when
$\sigma\in(1/2,1)$ or $\sigma=1/2$ and $\delta>1$, we can take
$\fl(t)\equiv 1$, so that for $\lambda$ large enough the solution of
(\ref{ODE-eqn})--(\ref{ODE-data}) is given by (\ref{formula:ul-sc}).
Thus, if we take $\tl:=(x_{2,\lambda})^{-1}$ and
$D:=(\delta+\sqrt{\delta^{2}-1})^{2}$, from (\ref{th:super-asympt}) we
obtain that
$$\lim_{\lambda\to
+\infty}\lambda\ul(\tl)=\left\{
\begin{array}{ll}
	1-e^{-1} & \mbox{if }\sigma>1/2  \\
	\noalign{\vspace{1ex}}
	1+(D-1)^{-1}(e^{-D}-De^{-1}) & \mbox{if $\sigma=1/2$ and $\delta>1$,}
\end{array}
\right.$$
and
$$\lim_{\lambda\to +\infty}\lambda^{\sigma}\ul'(\tl)=\left\{
\begin{array}{ll}
	(2\delta)^{-1} e^{-1} & \mbox{if }\sigma>1/2  \\
	\noalign{\vspace{1ex}}
	(2\sqrt{\delta^{2}-1})^{-1}(e^{-1}-e^{-D}) & 
	\mbox{if $\sigma=1/2$ and $\delta>1$.}
\end{array}
\right.$$

It is not difficult to see that these limits are always finite and
different from 0, which proves that $(\sigma,1,\sigma)$ is a blow up
triple in the case of supercritical dissipation with $\sigma<1$ (the
latter condition guarantees that $\tl\to 0^{+}$).

In the case of critical dissipation, namely when $\sigma=1/2$ and 
$\delta=1$, we can take once again $\fl(t)\equiv 1$. The 
solution of (\ref{ODE-eqn})--(\ref{ODE-data}) is
$$\ul(t)=\frac{1}{\lambda}-
\frac{1}{\lambda}\left(\lambda^{1/2}t+1\right)
e^{-\lambda^{1/2}t}.$$

Setting $\tl:=\lambda^{-1/2}$ we have that
$$\lim_{\lambda\to +\infty}\lambda\ul(\tl)=1-2e^{-1}
\quad\quad\mbox{and}\quad\quad
\lim_{\lambda\to +\infty}\lambda^{1/2}\ul'(\tl)=e^{-1},$$
which proves that $(1/2,1,1/2)$ is a blow-up triple in the critical 
case.

In the case of subcritical dissipation, namely when $\sigma<1/2$, or
$\sigma=1/2$ and $\delta\in(0,1)$, we need to use a family of
oscillating forcing terms.  Let us assume that $\lambda$ is large
enough so that the roots of the characteristic polynomial are of the
form $-\al\pm i\bl$, and let us set
\begin{equation}
	\tl:=\frac{W}{\al},
	\hspace{3em}
	\fl(t):=\sin(\bl(\tl-t)+\psi),
	\label{defn:tl-fl}
\end{equation}
where $W>0$ and $\psi\in\re$ are parameters to be chosen in the
sequel.  

As we already observed, when the roots of the characteristic
polynomial are complex conjugate numbers the solution $\ul(t)$ of
(\ref{ODE-eqn})--(\ref{ODE-data}) and its derivative $\ul'(t)$ are
given by (\ref{soln:ul}) and (\ref{soln:ul'}), respectively.
Keeping (\ref{defn:tl-fl}) into account, with the variable change 
$x=\al(\tl-s)$ we obtain that
\begin{equation}
	\lambda^{\sigma+1/2}\ul(\tl)=\frac{\lambda^{\sigma+1/2}}{\al\bl}
	\left(\cos\psi\cdot S_{\lambda}+
	\sin\psi\cdot M_{\lambda}\right),
	\label{lim:but-u}
\end{equation}
\begin{equation}
	\lambda^{\sigma}\ul'(\tl)=-\frac{\lambda^{\sigma}}{\bl}
	\left(\cos\psi\cdot S_{\lambda}+
	\sin\psi\cdot M_{\lambda}\right)+\frac{1}{\delta}
	\left(\cos\psi\cdot M_{\lambda}+\sin\psi\cdot C_{\lambda}\right),
	\label{lim:but-u'}
\end{equation}
where
$$S_{\lambda}:=\int_{0}^{W}e^{-x}\sin^{2}\left(\frac{\bl}{\al}x\right)dx,
\hspace{3em}
C_{\lambda}:=\int_{0}^{W}e^{-x}\cos^{2}\left(\frac{\bl}{\al}x\right)dx,$$
$$M_{\lambda}:=\int_{0}^{W}e^{-x}\sin\left(\frac{\bl}{\al}x\right)
\cos\left(\frac{\bl}{\al}x\right)dx.$$

Now we have to compute the limits as $\lambda\to +\infty$.  The limits
of the coefficients $\lambda^{\sigma+1/2}(\al\bl)^{-1}$ and
$\lambda^{\sigma}\bl^{-1}$ follow easily from the explicit
expressions (\ref{defn:al-bl-c}).  As for the integrals, we have to
distinguish two cases.

When $\sigma<1/2$, the coefficient $\bl\al^{-1}$ in the trigonometric
functions tends to $+\infty$.  This produces an homogenization
effect in the integrals, so that
$$\lim_{\lambda\to+\infty}S_{\lambda}=\lim_{\lambda\to+\infty}C_{\lambda}=
\frac{1}{2}\int_{0}^{W}e^{-x}\,dx=\frac{1}{2}\left(1-\frac{1}{e^{W}}\right),
\hspace{3em}
\lim_{\lambda\to+\infty}M_{\lambda}=0.$$

Plugging these limits into (\ref{lim:but-u}) and (\ref{lim:but-u'}),
we obtain that
$$\lim_{\lambda\to +\infty}\lambda^{\sigma+1/2}\ul(\tl)=
\frac{1}{2\delta}\left(1-\frac{1}{e^{W}}\right)\cos\psi,
\quad\quad
\lim_{\lambda\to +\infty}\lambda^{\sigma}\ul'(\tl)=
\frac{1}{2\delta}\left(1-\frac{1}{e^{W}}\right)\sin\psi.$$

Both limits are different from 0 for many values of $W$ and $\psi$
(for example $W=1$ and $\psi=\pi/4$), and this proves that
$(\sigma,\sigma+1/2,\sigma)$ is a blow-up triple for
$\sigma\in(0,1/2)$.

When $\sigma=1/2$, we have that both $\lambda^{\sigma+1/2}\ul(\tl)$ 
and $\lambda^{\sigma}\ul'(\tl)$ do not depend on $\lambda$. In 
particular, for $\psi=\pi/2$ we obtain that
$$\lambda^{\sigma+1/2}\ul(\tl)=\frac{1}{\delta\sqrt{1-\delta^{2}}}
\int_{0}^{W}e^{-x}\sin\left(Dx\right)
\cos\left(Dx\right)dx,$$
$$\lambda^{\sigma}\ul'(\tl)=\frac{1}{\delta}
\int_{0}^{W}e^{-x}\cos^{2}\left(Dx\right)dx-
\frac{1}{\sqrt{1-\delta^{2}}}
\int_{0}^{W}e^{-x}\sin\left(Dx\right)
\cos\left(Dx\right)dx,$$
where $D:=\delta^{-1}\sqrt{1-\delta^{2}}$.  Now it is easy to see that
both expressions are positive when $W>0$ is small enough.  This proves
that $(1/2,1,1/2)$ is a blow-up triple also in the subcritical case.

\paragraph{\textmd{\emph{Growth of one component in a short time}}}

Let $0\leq A<T$, let $(\sigma,\sigma_{0},\sigma_{1})$ be a triple of
positive real numbers satisfying the blow-up condition, and let
$c_{0}$ and $c_{1}$ be the two limits in (\ref{defn:buc-limits}).  Due
to the definition of blow-up triple, there exists $\Lambda\geq 0$ such
that the following
conditions 
$$\tl\leq T-A, 
\quad\quad\quad
\lambda^{\sigma_{0}}|\ul(\tl)|\geq\frac{3}{4}c_{0}, 
\quad\quad\quad
\lambda^{\sigma_{1}}|\ul'(\tl)|\geq\frac{3}{4}c_{1}$$
hold true for every $\lambda\geq\Lambda$. 

We claim that, for every $\lambda\geq\Lambda$, there exists
$B_{\lambda}\in(A,T)$ and a continuous function
$g_{\lambda}:[0,+\infty)\to\re$ such that
\begin{equation}
	|g_{\lambda}(t)|\leq 1
	\quad\quad
	\forall t\geq 0,
	\label{th:nr-T-1}
\end{equation}
\begin{equation}
	g_{\lambda}(t)=0
	\quad\quad
	\forall t\in[0,A]\cup[B_{\lambda},+\infty),
	\label{th:nr-T-2}
\end{equation}
and the unique solution $\ul(t)$ of (\ref{ODE-eqn})--(\ref{ODE-data})
with $g_{\lambda}(t)$ instead of $\fl(t)$ satisfies
\begin{equation}
	\lambda^{\sigma_{0}}|\ul(T)|\geq\frac{c_{0}}{2}
	\quad\quad\mbox{and}\quad\quad
	\lambda^{\sigma_{1}}|\ul'(T)|\geq\frac{c_{1}}{2}.
	\label{th:nr-T-3}
\end{equation}

In order to prove this result, we begin by defining a piecewise 
continuous function
$\phil:[0,+\infty)\to\re$ as
$$\phil(t):=\left\{
\begin{array}{ll}
	\fl(t-(T-\tl)) & \mbox{if }t\in[T-\tl,T],  \\
	\noalign{\vspace{0.5ex}}
	0 & \mbox{otherwise},
\end{array}
\right.$$
where $\fl(t)$ is the function which appears in the blow-up 
condition. It is easy to see that the solution of the ordinary 
differential equation
\begin{equation}
	\vl''(t)+2\delta \lambda^{\sigma}\vl'(t)+\lambda\vl(t)=\phil(t),	
	\label{pbm:eqn-v}
\end{equation}
with null initial data $\vl(0)=\vl'(0)=0$, is given by
$$\vl(t)=\left\{
\begin{array}{ll}
	0 & \mbox{if }t\in[0,T-\tl],  \\
	\noalign{\vspace{0.5ex}}
	\ul(t-(T-\tl)) & \mbox{if }t\in[T-\tl,T],
\end{array}
\right.$$
so that
\begin{equation}
	\lambda^{\sigma_{0}}|\vl(T)|=
	\lambda^{\sigma_{0}}|\ul(\tl)|\geq\frac{3}{4}c_{0}
	\quad\quad\mbox{and}\quad\quad
	\lambda^{\sigma_{1}}|\vl'(T)|=
	\lambda^{\sigma_{1}}|\ul'(\tl)|\geq\frac{3}{4}c_{1}.
	\label{est:c0-c1}
\end{equation}

Now we approximate $\phil(t)$ with suitable continuous functions 
$\phile(t)$. To this end, for every $\ep\in(0,\tl/4)$ we consider a 
cut-off function $\psie:[0,+\infty)\to[0,1]$ with 
$$\psie(t)=0
\quad\quad
\forall t\in [0,T-\tl+\ep]\cup[T-\ep,+\infty],$$
$$\psie(t)=1
\quad\quad
\forall t\in [T-\tl+2\ep,T-2\ep],$$
and then we set $\phile(t):=\psie(t)\cdot\phil(t)$. Let $\vle(t)$ be 
the solution of (\ref{pbm:eqn-v}), with $\phile(t)$ instead of 
$\phil(t)$, and null initial data.

It is easy to see that $\phile(t)\to\phil(t)$ in $L^{p}((0,+\infty))$
for every $p<+\infty$.  This is more than enough to guarantee that
$\vle\to\vl$ in the energy space, hence
$$\lim_{\ep\to 0^{+}}\vle(T)=\vl(T)
\quad\quad\mbox{and}\quad\quad
\lim_{\ep\to 0^{+}}\vle'(T)=\vl'(T).$$

Keeping (\ref{est:c0-c1}) into account, for $\ep(\lambda)$ small
enough we have that
$$
\lambda^{\sigma_{0}}|v_{\lambda,\ep(\lambda)}(T)|\geq\frac{c_{0}}{2}
\quad\quad\mbox{and}\quad\quad
\lambda^{\sigma_{1}}|v_{\lambda,\ep(\lambda)}'(T)|\geq\frac{c_{1}}{2}.
$$

Therefore, our requirements (\ref{th:nr-T-1}) through 
(\ref{th:nr-T-3}) are fulfilled by taking 
$g_{\lambda}(t):=\varphi_{\lambda,\ep(\lambda)}(t)$ and 
$B_{\lambda}:=T-\ep(\lambda)$.

\paragraph{\textmd{\emph{Lack of regularity for a given positive time}}}

Let $T>0$, let $\eta>0$, let $(\sigma,\sigma_{0},\sigma_{1})$ be a
triple of positive real numbers satisfying the blow-up condition, let
$\{\nu_{k}\}$ be any unbounded sequence of eigenvalues of $A$, let
$\{\widehat{e}_{k}\}$ be a corresponding sequence of orthonormal
eigenvectors, and let $\widehat{H}$ be the subspace of $H$ generated
by $\{\widehat{e}_{k}\}$.

We claim that there exists a function $f\in C^{0}_{b}([0,+\infty),H)$
satisfying (\ref{est:ft-eta}) and (\ref{est:ft-H}), and such that the
corresponding solution $u(t)$ of
(\ref{pbm:nh-eqn})--(\ref{pbm:nh-data}) satisfies
\begin{equation}
	u(T)\not\in D(A^{\sigma_{0}})
	\quad \mbox{and}\quad 
	u'(T)\not\in D(A^{\sigma_{1}}).
	\label{not-sigma-01}
\end{equation}

In order to prove this claim, we begin by choosing a sequence 
$\{\omega_{n}\}\subseteq[0,1]$ of positive real numbers such that
$$\lim_{n\to +\infty}\omega_{n}=0
\quad\quad\mbox{and}\quad\quad
\sum_{n=0}^{\infty}\omega_{n}^{2}=+\infty.$$

Then we choose an increasing sequence $\{k_{n}\}$ of positive
integers, an increasing sequence $\{B_{n}\}\subseteq[0,T)$ of times,
and a sequence $f_{n}:[0,+\infty)\to[0,1]$ of continuous functions
such that $f_{n}(t)=0$ for every $t\not\in(B_{n-1},B_{n})$ (hence with
disjoint supports), and such that the corresponding solutions $u_{n}(t)$
of (\ref{pbm:nh-eqn})--(\ref{pbm:nh-data}) satisfy
$$
\nu_{k_{n}}^{\sigma_{0}}|u_{n}(T)|\geq\frac{c_{0}}{2}
\quad\quad\mbox{and}\quad\quad
\nu_{k_{n}}^{\sigma_{1}}|u_{n}'(T)|\geq\frac{c_{1}}{2}.
$$

As soon as we have such sequences, our claim follows with
$$f(t):=\eta\sum_{n=1}^{\infty}\omega_{n}f_{n}(t)\widehat{e}_{k_{n}}.$$

Indeed the series converges because its terms have disjoint supports 
and their norm goes to 0. The sum $f(t)$ satisfies (\ref{est:ft-H}) 
for trivial reasons, and satisfies (\ref{est:ft-eta}) because all 
terms do and have disjoint supports. The corresponding solution $u(t)$ 
of (\ref{pbm:nh-eqn})--(\ref{pbm:nh-data}) is clearly
$$u(t):=\eta\sum_{n=1}^{\infty}\omega_{n}u_{n}(t)\widehat{e}_{k_{n}},$$
so that
$$\left|A^{\sigma_{0}}u(T)\right|^{2}=
\eta^{2}\sum_{n=1}^{\infty}\omega_{n}^{2}\nu_{k_{n}}^{2\sigma_{0}}|u_{n}(T)|^{2}
\geq\eta^{2}\frac{c_{0}^{2}}{4}\sum_{n=1}^{\infty}\omega_{n}^{2}
=+\infty,$$
$$\left|A^{\sigma_{1}}u'(T)\right|^{2}=
\eta^{2}\sum_{n=1}^{\infty}\omega_{n}^{2}\nu_{k_{n}}^{2\sigma_{1}}|u_{n}'(T)|^{2}
\geq\eta^{2}\frac{c_{1}^{2}}{4}\sum_{n=1}^{\infty}\omega_{n}^{2}
=+\infty,$$
which proves (\ref{not-sigma-01}).

In order to define the sequences we need, we repeatedly apply the
result of the previous paragraph.  First of all, we apply it with
$A:=0$ and we obtain a value $\Lambda_{1}$ such that for every
$\lambda\geq\Lambda_{1}$ there exists $B_{\lambda}\in(0,T)$ and
$g_{\lambda}:[0,+\infty)\to[0,1]$ satisfying (\ref{th:nr-T-1}) through
(\ref{th:nr-T-3}).  Now we choose a positive integer $k_{1}$ such that
$\nu_{k_{1}}\geq\Lambda_{1}$, and then we set
$B_{1}:=B_{\nu_{k_{1}}}$ and $f_{1}(t):=g_{\nu_{k_{1}}}(t)$.

Then we proceed by induction.  Let us assume that $k_{n-1}$, $B_{n-1}$
and $f_{n-1}(t)$ have been defined, and let us apply the result of the
previous paragraph with $A:=B_{n-1}$.  Once again we obtain a value
$\Lambda_{n}$ such that for $\lambda\geq\Lambda_{n}$ there exists
$B_{\lambda}\in(B_{n-1},T)$ and $g_{\lambda}:[0,+\infty)\to[0,1]$
satisfying (\ref{th:nr-T-1}) through (\ref{th:nr-T-3}).  Now we choose
a positive integer $k_{n}>k_{n-1}$ such that
$\nu_{k_{n}}\geq\Lambda_{n}$, and then we set $B_{n}:=B_{\nu_{k_{n}}}$
and $f_{n}(t):=g_{\nu_{k_{n}}}(t)$.

\paragraph{\textmd{\emph{Lack of regularity for a given sequence of times}}}

In the previous paragraph we produced a solution of
(\ref{pbm:nh-eqn})--(\ref{pbm:nh-data}) which has the required
regularity loss at a given time $T>0$.  This solution has been
constructed using a given unbounded sequence of eigenvalues of $A$,
and an external force as small as we want.  Now we need to produce the
same regularity loss for all times in a given sequence $\{t_{n}\}$.

The procedure is exactly the same as in the last paragraph of
Section~\ref{ctrex:stat-1}.  We partition the sequence
$\{\lambda_{k}\}$ into countably many disjoint subsequences, and for
each $n\in\n$ we apply the previous result in order to obtain an
external force $f_{n}(t)$, with norm less than or equal to $2^{-n}$,
for which the corresponding solution $u_{n}(t)$ of
(\ref{pbm:nh-eqn})--(\ref{pbm:nh-data}) has the required regularity
loss at time $t_{n}$.  The series of these external forces converges
to an external force with all the properties we need.  We refer to the
last paragraph of Section~\ref{ctrex:stat-1} for the details.

\subsection{Proof of~(\ref{th:4b}) in statement~(3)}

The construction is analogous to the one for statement~(2), the only
difference being that now we have to produce a regularity loss only
for the derivative.

We begin by saying that a pair $(\sigma,\sigma_{1})$ of positive real
numbers satisfies the blow-up condition for the derivative if there
exist families $\{\tl\}\subseteq(0,+\infty)$ and $\{\fl\}\subseteq
C^{0}_{b}([0,+\infty),H)$ satisfying (\ref{defn:buc-fl}) and
(\ref{defn:buc-tl}), and such that the corresponding solutions
$\ul(t)$ of (\ref{ODE-eqn})--(\ref{ODE-data}) satisfy the second
relation in~(\ref{defn:buc-limits}).

Then we prove that the pairs $\left(\sigma,\sigma\right)$ satisfy the
blow-up condition for the derivative for every $\sigma\geq 1$.
Indeed, if we take $\fl(t)\equiv 1$, for $\lambda$ large enough the
solution of (\ref{ODE-eqn})--(\ref{ODE-data}) is given by
(\ref{formula:ul-sc}), hence its derivative is
$$\ul'(t)=\frac{1}{x_{1,\lambda}-x_{2,\lambda}}
\left(e^{-x_{2,\lambda}t}-e^{-x_{1,\lambda}t}\right).$$

At this point the conclusion easily follows with
$\tl:=(x_{1,\lambda})^{-1}$.

From now on the argument is exactly the same as in the proof of
statement~(2).  First we obtain the growth of (the derivative of) one
component in an arbitrary short time, then the lack of regularity at
a given positive time, and finally the lack of regularity at a given
sequence of positive times.

\subsection{Proof of statement~(4)}

\paragraph{\textmd{\textit{Preliminary integral estimate}}}

Let $\{\alpha_{n}\}$ be any sequence of positive real numbers such that 
$\alpha_{n}\to 0^{+}$. We claim that there exist an increasing 
sequence $\{k_{n}\}$ of positive integers, and an increasing sequence 
$\{T_{n}\}\subseteq[0,+\infty)$ with $T_{0}=0$ such that
\begin{equation}
	\alpha_{k_{n}}\int_{T_{n-1}}^{T_{n}}e^{-\alpha_{k_{n}}x}dx\geq
	\frac{1}{e}\left(1-\frac{1}{e}\right)
	\quad\quad
	\forall n\geq 1.
	\label{th:alpha-T}
\end{equation}

In order to prove this claim, we set $k_{1}:=1$, and then by induction we choose 
$k_{n+1}>k_{n}$ in such a way that
\begin{equation}
	\frac{1}{\alpha_{k_{n+1}}}\geq
	\frac{1}{\alpha_{k_{1}}}+\ldots+\frac{1}{\alpha_{k_{n}}}
	\quad\quad
	\forall n\geq 1.
	\label{defn:kn}
\end{equation}

Such a choice is possible because $\alpha_{n}\to 0^{+}$. Then we set 
$T_{0}:=0$ and
\begin{equation}
	T_{n}:=\frac{1}{\alpha_{k_{1}}}+\ldots+\frac{1}{\alpha_{k_{n}}}
	\quad\quad
	\forall n\geq 1.
	\label{defn:Tn}
\end{equation}

Due to (\ref{defn:kn}) and (\ref{defn:Tn}) we have that 
$\alpha_{k_{n}}T_{n-1}\leq 1$ and $\alpha_{k_{n}}(T_{n}-T_{n-1})=1$, hence
$$\alpha_{k_{n}}\int_{T_{n-1}}^{T_{n}}e^{-\alpha_{k_{n}}x}dx=
e^{-\alpha_{k_{n}}T_{n-1}}\left(1-
e^{-\alpha_{k_{n}}(T_{n}-T_{n-1})}\right)\geq
\frac{1}{e}\left(1-\frac{1}{e}\right)$$
for every $n\geq 1$, as required.

\paragraph{\textmd{\textit{Passing any given threshold}}}

Let $M\geq 0$, let $\eta>0$, let $\{\nu_{n}\}$ be any unbounded
sequence of eigenvalues of $A$, and let $\{\widehat{e}_{n}\}$ be a corresponding
sequence of orthonormal eigenvectors.  

We claim that there exist $T>0$, a subspace $\widehat{H}$ of $H$
generated by a finite subset of $\{\widehat{e}_{n}\}$, and a function $f\in
C^{0}_{b}([0,+\infty),H)$ satisfying (\ref{est:ft-eta}) and
(\ref{est:ft-H}), and such that the corresponding solution $u(t)$ of
(\ref{pbm:nh-eqn})--(\ref{pbm:nh-data}) satisfies
\begin{equation}
	|Au(T)|^{2}\geq M.
	\label{th:Au(T)}
\end{equation}

In other words, the external force is as small as we want and
concentrated on a finite number of components, but $|Au(T)|$ exceeds a
given threshold.

In order to prove the claim, let us consider the roots of the
characteristic polynomial~(\ref{char-pol}) with $\lambda=\nu_{n}$, and
let us set for simplicity $x_{1,n}:=x_{1,\nu_{n}}$ and
$x_{2,n}:=x_{2,\nu_{n}}$.  We always assume that $\nu_{n}$ is large
enough so that these roots are distinct real numbers.  We
also assume also that $\nu_{n}$ is large enough so that
\begin{equation}
	\frac{\nu_{n}}{x_{1,n}-x_{2,n}}\leq 1,
	\quad\quad\quad
	x_{1,n}\geq 1,
	\quad\quad\quad
	\frac{\nu_{n}}{x_{1,n}-x_{2,n}}\cdot\frac{1}{x_{2,n}}\geq\frac{1}{2}.
	\label{hp:nu-n}
\end{equation}

This is clearly possible because of (\ref{th:super-asympt}).  Since
$\sigma>1$, we have also that $x_{2,n}\to 0^{+}$, and therefore we can
apply the result of the previous paragraph with $\alpha_{n}:=x_{2,n}$.
Let $k_{n}$ and $T_{n}$ be the corresponding sequences for which
(\ref{th:alpha-T}) holds true.

Let us consider the piecewise constant function
$\psi:[0,+\infty)\to\re$ defined by $\psi(t):=\widehat{e}_{k_{n}}$ for every
$t\in[T_{n-1},T_{n})$ and every $n\in\n$.  Let us choose a positive
integer $N$ large enough so that
$$\eta^{2}\frac{1}{4e^{2}}\left(1-\frac{1}{e}\right)^{2}\cdot N
\geq (2M+\eta)^{2}.$$

Let us set $T:=T_{N}$ and
$$g(t):=\left\{
\begin{array}{ll}
	\eta\psi(T-t) & \mbox{if }t\in[0,T],  \\
	\noalign{\vspace{0.5ex}}
	0 & \mbox{if }t>T.
\end{array}
\right.$$

Let $\widehat{H}$ be the subspace of $H$ generated by
$\{\widehat{e}_{k_{1}},\ldots,\widehat{e}_{k_{N}}\}$.  The function $g(t)$ is not
continuous, but it satisfies (\ref{est:ft-eta}) and (\ref{est:ft-H}).
The corresponding solution of (\ref{pbm:nh-eqn})--(\ref{pbm:nh-data})
is
$$u(t):=v(t)+w(t):=\eta\sum_{n=1}^{N}v_{n}(t)\widehat{e}_{k_{n}}+
\eta\sum_{n=1}^{N}w_{n}(t)\widehat{e}_{k_{n}},$$ 
with
$$v_{n}(t):=-\frac{1}{x_{1,k_{n}}-x_{2,k_{n}}}\int_{0}^{t}
e^{-x_{1,k_{n}}(t-s)}\psi_{k_{n}}(T-s)\,ds,$$
$$w_{n}(t):=\frac{1}{x_{1,k_{n}}-x_{2,k_{n}}}\int_{0}^{t}
e^{-x_{2,k_{n}}(t-s)}\psi_{k_{n}}(T-s)\,ds,$$
where of course $\psi_{k_{n}}(t)$ denotes the component of $\psi(t)$ with
respect to $\widehat{e}_{k_{n}}$.

Let us estimate $v(T)$ and $w(T)$ separately. 
In order to estimate $v(T)$, we argue as in the proof of
Lemma~\ref{lemma:main}.  We consider $T$ as a parameter and we
introduce the vector
$$\varphi(s):=\sum_{n=1}^{N}\varphi_{n}(s)\widehat{e}_{k_{n}}$$ 
with components
$$\varphi_{n}(s):=-\frac{\nu_{k_{n}}}{x_{1,k_{n}}-x_{2,k_{n}}}
e^{-x_{1,k_{n}}(T-s)}\psi_{k_{n}}(T-s),$$
so that
$$|Av(T)|=\eta\left|\int_{0}^{T}\varphi(s)\,ds\right|\leq
\eta\int_{0}^{T}\left|\varphi(s)\right|\,ds.$$

In order to estimate $|\varphi(s)|$, we exploit the first two
conditions in (\ref{hp:nu-n}) and we obtain that $|\varphi_{n}(s)|\leq
e^{-(T-s)}$ for every $n=1,\ldots,N$ and every $s\geq 0$.  Since only
one component of $\varphi(s)$ is different from 0 for each $s$, we
conclude that $|\varphi(s)|\leq e^{-(T-s)}$, hence
\begin{equation}
	|Av(T)|\leq\eta\int_{0}^{T}e^{-(T-s)}\,ds\leq\eta.
	\label{est:v(T)}
\end{equation}

In order to estimate $w(T)$, we first observe that
$$\int_{0}^{T}e^{-x_{2,k_{n}}(T-s)}\psi_{k_{n}}(T-s)\,ds=
\int_{0}^{T}e^{-x_{2,k_{n}}y}\psi_{k_{n}}(y)\,dy=
\int_{T_{n-1}}^{T_{n}}e^{-x_{2,k_{n}}y}\,dy.$$

From (\ref{th:alpha-T}) and the last condition in (\ref{hp:nu-n}) we 
deduce that
\begin{eqnarray*}
	\nu_{k_{n}}w_{n}(T) & = & \frac{\nu_{k_{n}}}{x_{1,k_{n}}-x_{2,k_{n}}}
	\int_{T_{n-1}}^{T_{n}}e^{-x_{2,k_{n}}y}\psi_{k_{n}}(y)\,dy \\
	 & \geq &
	 \frac{\nu_{k_{n}}}{x_{1,k_{n}}-x_{2,k_{n}}}\cdot\frac{1}{x_{2,k_{n}}}
	\cdot\frac{1}{e}\left(1-\frac{1}{e}\right) \\
	 & \geq & \frac{1}{2e}\left(1-\frac{1}{e}\right),
\end{eqnarray*}
and finally
\begin{equation}
	|Aw(T)|^{2}=\eta^{2}\sum_{n=1}^{N}|w_{n}(T)|^{2}\geq
	\eta^{2}\frac{1}{4e^{2}}\left(1-\frac{1}{e}\right)^{2}\cdot N
	\geq(2M+\eta)^{2}.
	\label{est:w(T)}
\end{equation}

Therefore, from (\ref{est:v(T)}) and (\ref{est:w(T)}) we conclude that
$$|Au(T)|\geq|Aw(T)|-|Av(T)|\geq 2M.$$

It remains to fix the issue that $g(t)$ is not continuous.  To this
end, it is enough to approximate $g(t)$ with a continuous function
$f(t)$ which still satisfies (\ref{est:ft-eta}) and (\ref{est:ft-H})
(to this end, it is enough to approximate from below the
characteristic functions of the intervals in the definition of
$\psi(t)$).  If $f(t)$ is close enough to $g(t)$, for example in
$L^{2}((0,T),H)$, then the corresponding solution of
(\ref{pbm:nh-eqn})--(\ref{pbm:nh-data}) is as close as we want to
$u(t)$ in the energy norm, which is equivalent to the norm in $D(A)$
because only a finite number of components is involved.  This proves
that we can choose $f(t)$ so that the new solution satisfies
(\ref{th:Au(T)}).

\paragraph{\textmd{\textit{Conclusion}}} 

We construct a sequence $\{t_{n}\}$ of positive times, a sequence
$\{H_{n}\}$ of (finite dimensional) pairwise orthogonal subspaces of
$H$, and a sequence of continuous functions $f_{n}:[0,+\infty)\to H$
such that $|f_{n}(t)|\leq 2^{-n}$ and $f_{n}(t)\in H_{n}$ for every
$t\geq 0$ and $n\in\n$, and such that the corresponding solutions
$u_{n}(t)$ of (\ref{pbm:nh-eqn})--(\ref{pbm:nh-data}) satisfy
$$|Au_{n}(t_{n})|^{2}\geq n
\quad\quad
\forall n\in\n.$$

The existence of such sequences follows easily from a repeated
application of the result of the previous paragraph, each time with
$\eta:=2^{-n}$, $M:=n$, and $\{\nu_{n}\}$ equal to the elements of the
sequence $\{\lambda_{n}\}$ which have not yet been used up to that
point (we recall that at each step only a finite number of eigenvalues
is involved).

The conclusion follows as in the previous cases by defining $f(t)$ and
$u(t)$ as in~(\ref{defn:f-u}).  Since the subspaces $H_{n}$ are
pairwise orthogonal, we have that $|Au(t_{n})|\geq|Au_{n}(t_{n})|\geq
n$ for every $n\in\n$, which proves (\ref{th:5}).\qed

\subsubsection*{\centering Acknowledgments}

This work hes been done while the first two authors where visiting the
Laboratoire Jacques Louis Lions of the UPMC (Paris~VI). The stay
was partially supported by the FSMP (Fondation Sciences
Math\'{e}matiques de Paris).

\label{NumeroPagine}


\begin{thebibliography}{99}
	
 	\bibitem{AP} \textsc{L.~Amerio, G.~Prouse}; Uniqueness and
	almost-periodicity theorems for a non linear wave equation.
	\emph{Atti Accad.\ Naz.\ Lincei Rend.\ Cl.\ Sci.\ Fis.\ Mat.\
	Natur.\ (8)} \textbf{46} (1969), 1--8.

	\bibitem{Bi} \textsc{M.~Biroli}; Bounded or almost periodic
	solution of the non linear vibrating membrane equation.
	\emph{Ricerche Mat.}\ \textbf{22} (1973), 190--202.

	\bibitem{B-H} \textsc{M.~Biroli, A.~Haraux}; Asymptotic behavior
	for an almost periodic, strongly dissipative wave equation.
	\emph{J.\ Differential Equations} \textbf{38} (1980), no.~3, 422--440.
	
	\bibitem {C-H}{\sc T.~Cazenave, A.~Haraux}; An introduction to
	semilinear evolution equations, \emph{Oxford Lecture Series in
	Mathematics and its Applications} \textbf{13}. The Clarendon Press,
	Oxford University Press, New York, 1998.

	\bibitem{CR}\textsc{G.~Chen, D.~L.~Russell}; A mathematical model
	for linear elastic systems with structural damping.  \emph{Quart.\
	Appl.\ Math.}\ \textbf{39} (1981/82), no.~4, 433--454.

	\bibitem{CT1}\textsc{S.~P.~Chen, R.~Triggiani}; Proof of
	extensions of two conjectures on structural damping for elastic
	systems.  \emph{Pacific J.\ Math.}\ \textbf{136} (1989), no.~1,
	15--55.

	\bibitem{CT2}\textsc{S.~P.~Chen, R.~Triggiani}; Characterization
	of domains of fractional powers of certain operators arising in
	elastic systems, and applications.  \emph{J.\ Differential
	Equations} \textbf{88} (1990), no.~2, 279--293.
	
	\bibitem{CT3}\textsc{S.~P.~Chen, R.~Triggiani}; Gevrey class
	semigroups arising from elastic systems with gentle dissipation:
	the case $0<\alpha<1/2$.  \emph{Proc.\ Amer.\ Math.\ Soc.}\
	\textbf{110} (1990), no.~2, 401--415.
	
	\bibitem{DAR}\textsc{M.\ D'Abbicco, M.\ Reissig}; Semi-linear
	structural damped waves.  \emph{Math.\ Methods Appl.\ Sci.}\ To
	appear.  (doi: 10.1002/mma.2913, \texttt{arXiv:1209.3204}).
	
	\bibitem{FGR}\textsc{L.~H.~Fatori, M.~Z.~Garay, J.~E.~M.~Rivera};
	Differentiability, analyticity and optimal rates of decay for
	damped wave equations.  \emph{Electron.\ J.\ Differential
	Equations} \textbf{2012}, No.~48, 13 pp.

	\bibitem{h-brezil}\textsc{A.~Haraux}; Uniform decay and Lagrange
	stability for linear contraction semi-groups.  \emph{Mat.\ Apl.\
	Comput.}\ \textbf{7} (1988), no.~3, 143--154.

	\bibitem{HO}\textsc{A.~Haraux, M.~\^{O}tani}; Analyticity and
	regularity for a class of second order evolution equation.
	\emph{Evol.\ Equat.\ Contr.\ Theor.}\ \textbf{2} (2013),
	no.~1, 101--117.
	
	\bibitem{H0}\textsc{A.~Haraux}; Nonlinear evolution
	equations--global behavior of solutions.  \emph{Lecture Notes in
	Mathematics} \textbf{841}.  Springer-Verlag, Berlin-New York, 1981.

	\bibitem{H4}\textsc{A.~Haraux}; Damping out of transient states
	for some semilinear, quasiautonomous systems of hyperbolic type.
	\emph{Rend.\ Accad.\ Naz.\ Sci.\ XL Mem.\ Mat.\ (5)} \textbf{7} (1983),
	89--136.

	\bibitem{H2}\textsc{A.~Haraux}; Semi-linear hyperbolic problems in
	bounded domains.  \emph{Math.\ Rep.}\ \textbf {3} (1987), no.~1, i--xxiv
	and 1--281.

	\bibitem{Anti-P}\textsc{A.~Haraux}; Anti-periodic solutions of
	some nonlinear evolution equations.  \emph{Manuscripta Math.}\
	\textbf{63} (1989), no.~4, 479--505.

	\bibitem{Reg-Forc}\textsc{A.~Haraux}; Nonresonance for a strongly
	dissipative wave equation in higher dimensions.  \emph{Manuscripta
	Math.}\ \textbf{53} (1985), no.~1-2, 145--166.

	\bibitem{HZ} \textsc{A.~Haraux, E.~Zuazua}; Decay estimates for
	some semilinear damped hyperbolic problems.  \emph{Arch.\ Rational
	Mech.\ Anal.}\ \textbf{100} (1988), no.~2, 191--206.

	\bibitem{I1}\textsc{R.~Ikehata}; Decay estimates of solutions for
	the wave equations with strong damping terms in unbounded domains.
	\emph{Math.\ Methods Appl.\ Sci.}\ \textbf{24} (2001), no.~9,
	659--670.

	\bibitem{I2}\textsc{R.~Ikehata, M.~Natsume}; Energy decay estimates
	for wave equations with a fractional damping.  \emph{ Differential
	Integral Equations} \textbf{25} (2012), no.~9-10, 939--956.

	\bibitem{I3}\textsc{R.~Ikehata, G.~Todorova, B.~Yordanov}; Wave
	equations with strong damping in Hilbert spaces.  \emph{J.\
	Differential Equations} \textbf{254} (2013), no.~8, 3352--3368.
	
	\bibitem{LZ}\textsc{B.~M.~Levitan, V.~V.~Zhikov}; Almost periodic
	functions and differential equations.  Translated from the Russian
	by L.~W.~Longdon.  \emph{Cambridge University Press, Cambridge-New
	York, 1982.  }

	\bibitem{LL}\textsc{K.~Liu, Z.~Liu}; Analyticity and
	differentiability of semigroups associated with elastic systems
	with damping and gyroscopic forces.  \emph{J.\ Differential
	Equations} \textbf{141} (1997), no.~2, 340--355.

	\bibitem{MR:gevrey}\textsc{S.\ Matthes, M.\ Reissig};
	Qualitative properties of structural damped wave models.
	\emph{Eurasian Math.\ J.}\ To appear.

	\bibitem{mugnolo}\textsc{D.~Mugnolo}; A variational approach to
	strongly damped wave equations.  \emph{Functional analysis and
	evolution equations}, 503--514, Birkh\"{a}user, Basel, 2008.

	\bibitem{nishihara}\textsc{K.~Nishihara}; Degenerate quasilinear
	hyperbolic equation with strong damping.  \emph{Funkcial.\
	Ekvac.}\ \textbf{27} (1984), no.~1, 125--145.

	\bibitem{ono-nishihara}\textsc{K.~Ono, K.~Nishihara}; On a
	nonlinear degenerate integro-differential equation of hyperbolic
	type with a strong dissipation.  \emph{Adv.\ Math.\ Sci.\ Appl.}\
	\textbf{5} (1995), no.~2, 457--476.
	
	\bibitem{reed}{\sc M.\ Reed, B.\ Simon}; \emph{Methods of Modern
	Mathematical Physics, I: Functional Analysis.  Second edition}.
	Academic Press, New York, 1980.
	
	\bibitem{shibata}\textsc{Y.~Shibata}; On the rate of decay of
	solutions to linear viscoelastic equation.  \emph{Math.\ Methods
	Appl.\ Sci.}\ \textbf{23} (2000), no.~3, 203--226.
	
\end{thebibliography}
\end{document}